\documentclass[12pt]{amsart}
\usepackage{amssymb}
\usepackage{graphics}

\begin{document}

\newtheorem{theorem}{Theorem}[section]
\newtheorem{lemma}[theorem]{Lemma}
\newtheorem{corollary}[theorem]{Corollary}
\newtheorem{proposition}[theorem]{Proposition}

\newtheorem{maintheorem}{Main Theorem}
\def\themaintheorem{\unskip}

\theoremstyle{definition}
\newtheorem{remark}{Remark}
\def\theremark{\unskip}
\newtheorem{definition}{Definition}
\def\thedefinition{\unskip}
\newtheorem{problem}{Problem}[section]

\numberwithin{equation}{section}

\def\Re{\operatorname{Re\,} }
\def\Im{\operatorname{Im\,} }
\def\dist{\operatorname{dist\,} }
\def\p{\partial}
\def\rp{^{-1}}
\def\esslimsup{\operatorname{esslimsup}}
\def\essliminf{\operatorname{essliminf}}

\def\iwithdot{i}
\def\jwithdot{j}

\def\gap{\smallskip\noindent}
\def\mgap{\medskip\noindent}
\def\hils{{\mathcal H}}
\def\R{{\Bbb R}}
\def\Z{{\Bbb Z}}
\def\Q{{\Bbb Q}}
\def\C{{\Bbb C}}
\def\knj{k^{(n)}_j}
\def\kmj{k^{(m)}_j}
\def\tnj{\theta^{(n)}_j}
\def\tmj{\theta^{(m)}_j}

\def\reals{ {{\mathbb R}} }
\def\naturals{ {{\mathbb N}} }
\def\integers{ {{\mathbb Z}} }
\def\complex{{\mathbb C}}
\newcommand\beq{\begin{equation}}
\newcommand\eeq{\end{equation}}
\def\scriptd{ {\mathcal D} }
\def\scripte{ {\mathcal E} }
\def\scriptb{ {\mathcal B} }
\def\scripth{ {\mathcal H} }
\def\scriptdplus{ {\mathcal D}_+ }
\def\scriptdminus{ {\mathcal D}_- }
\def\scriptdpm{ {\mathcal D}_\pm }
\def\scriptdmp{ {\mathcal D}_\mp }
\def\scriptm{{\mathcal M}}
\def\phila{\phi_\lambda}
\def\scriptt{{\mathcal T}}
\def\scriptf{{\mathcal F}}
\def\scriptg{{\mathcal G}}
\def\scriptv{{\mathcal V}}
\def\dthk{{\frac{\partial \theta}{\partial k}}}
\def\sump{\sideset{}{'}\sum}
\def\marrow{\overset{\rightharpoonup}{m}}
\def\jarrow{\overset{\rightharpoonup}{\jwithdot}}

\author{Alexander Kiselev}
\address{ Department of Mathematics \\
University of Chicago \\
Chicago, IL 60637}
\email{kiselev@math.uchicago.edu}
\thanks{The author was supported in part by NSF grant DMS-0102554 and 
by Alfred P. Sloan Research Fellowship}



\title[Singular Continuous Spectrum]
{Imbedded Singular Continuous Spectrum for Schr\"odinger Operators}

\begin{abstract}
We construct examples of potentials $V(x)$ satisfying $|V(x)| \leq \frac{h(x)}{1+x},$ 
where the function $h(x)$ is growing arbitrarily slowly, such that the corresponding Schr\"odinger
operator has imbedded singular continuous spectrum. This solves one of the fifteen ``twenty-first century" 
problems for Schr\"odinger operators  posed 
by Barry Simon in \cite{BS}. The construction also provides the first example of a Schr\"odinger operator for
which  M\"oller wave operators exist but are not  asymptotically complete due to the presence of singular continuous spectrum.
\end{abstract}

\maketitle

\section{Introduction} \label{intro}

Let $H_V$ denote one-dimensional Schr\"odinger operator defined on the half-axis 
by differential expression
\begin{equation}\label{so}
 H_V = -\frac{d^2}{dx^2} +V(x) 
\end{equation}
and some self-adjoint boundary condition at  zero.  
The operator \eqref{so} describes a charged particle, such as an electron, in electric 
field $V(x).$ When $V(x)$ is decaying quickly,
one expects the spectral and dynamical properties of $H_V$ to remain close to 
those of free operator $H_0.$ Recall that with the operator \eqref{so} one can associate
in a canonical way a spectral measure $\mu$ (see, for example, \cite{CL,CarLa}) which contains 
much information about the quantum system. 
The classical result in this direction, going back to the beginning of the century,
 is that if $V \in L^1,$ the spectral measure on the positive semi-axis
remains purely absolutely continuous.
In one dimension, absolutely continuous spectrum corresponds to the ballistic rate of propagation of the particle.  
A natural question to ask is what are the critical rates of decay at which some changes in the 
spectral and dynamical properties of $H_V$ may happen. As far as $V(x)$ has any decay at all, 
the essential spectrum coincides with $[0,\infty),$ but the quality of the spectrum and dynamics may change.
In 1928, Wigner and von Neumann \cite{WvN} showed that there exist potentials $V(x)$ 
satisfying $|V(x)| \leq \frac{C}{x+1}$ such that $H_V$ has positive imbedded eigenvalue $E=1$. 
This is a purely quantum resonance phenomenon, since $V$ can be chosen arbitrarily small and bound state 
originates from long range correlations in the potential rather than usual confining effect. 
Naboko \cite{Nab} and Simon \cite{Sim1} provided constructions which show that much more 
drastic changes are possible if potential decays arbitrarily slower than Coulomb. Namely, for any 
positive function $h(x)$ which grows at infinity, there exist potentials $V$ such that $|V(x)| 
\leq \frac{h(x)}{1+x}$ and $H_V$ has a dense set of positive eigenvalues.  
Due to some recent results \cite{CK1,Rem1}, it is known that the absolutely continuous 
spectrum on the positive semi-axis is preserved for $|V(x)| \leq C(1+x)^{-1/2-\epsilon}$ and 
in fact for $V \in L^2$ \cite{DK}. More precisely, for such potentials the absolutely continuous part of the 
spectral measure, $\mu^{{\rm ac}},$ gives positive weight to any subset of $(0,\infty)$ of 
positive Lebesgue measure. 
Therefore, Naboko and Simon examples  have dense point spectrum 
imbedded into the absolutely continuous spectrum. 

In \cite{BS}, Simon has posed fifteen problems on Schr\"odinger operators, one of which asks whether 
Schr\"odinger operators with potentials satisfying $|V(x)| \leq C(1+|x|)^{-1/2-\epsilon}$ can also 
have imbedded singular continuous spectrum. 
Remling \cite{Re2} and later Christ-Kiselev \cite{CK2} provided estimates on the dimension of the support of 
the possible singular component of the measure for $|V(x)| \leq C(1+x)^{-\alpha}.$ Remling \cite{Re3} and 
Kriecherbauer and Remling \cite{KR1} 
constructed power-decaying potentials which have decaying solutions on the sets of Hausdorff 
dimension exactly as predicted by \cite{Re2,CK2}.  However these constructions are not sufficient to infer
the  presence of the singular continuous component. The problem of controlling singular continuous spectrum 
is difficult since one has to study a set of energies which is uncountable but has Lebesgue measure zero. 
Because of that, in most known situations, the 
existence of singular continuous spectrum is proven ``not by what it is, but by what it is not":
by elimination of possibilities to have absolutely continuous or point spectrum. 
For instance, see the classical example of Pearson \cite{Pearson} or genericity of singular continuous spectrum results 
of Simon \cite{Sim3}. 
However, when singular 
continuous spectrum is imbedded, such path is obviously ruled out. 
 To the best of our knowledge, up until now there were no constructive examples 
of Schr\"odinger operators with imbedded singular continuous spectrum in any setting. 
We notice that the fact that potentials leading to such operators exist  
follows from the inverse spectral theory (see, e.g. \cite{Lev,Mar}). This classical result 
has been recently improved by Denisov \cite{Den} to give existence of potentials leading to imbedded singular
continuous spectrum in the $L^2$ class. Killip and Simon \cite{KilSi} have subsequently found necessary and sufficient 
conditions for a measure to be a spectral measure of a free Jacobi matrix with $l^2$ perturbation. However these
interesting results provide only indirect information on what a potential leading to singular continuous spectrum looks like
and do not seem promising in terms of obtaining results for power decaying potentials.

There is an additional point of view on the problem of imbedded singular continuous spectrum for potentials which 
decay faster than $L^2$ which in our opinion makes the question even more interesting.   
Recall the definition of wave operators 
\begin{equation}\label{wo}
\Omega_\pm f = s-\lim_{t \mp \infty} e^{iH_V t} e^{-iH_0 t} f, 
\end{equation}
where the limit is taken in the strong $L^2$ sense. Existence of wave operators implies that the absolutely 
continuous spectrum of $H_V$ fills the whole positive semi-axis, and provides much more precise information about 
long-time perturbed dynamics. The wave operators are called asymptotically complete (see, e.g. \cite{RS3}) 
if the range of $\Omega_\pm$ 
coincides with the orthogonal complement of the subspace spanned by eigenvectors of
the operator $H_V$. An alternative equivalent characterization is that the range of 
the wave operators
is equal to the absolutely continuous subspace $\scripth_{\text{ac}}(H_V)$ of the operator
$H_V$, and that the singular continuous spectrum $\sigma_{{\rm sc}}(H_V)$ is empty.  
The intuitive meaning of asymptotic completeness 
is that the dynamics 
of the perturbed operator can be divided into two well-understood parts: 
scattering states traveling 
to infinity in a way similar to the free evolution, and bound states which
 remain confined in a certain sense
for all times.
On the other hand, the singular continuous spectrum is less understood, and corresponds
to an intermediate situation: states which do not remain bounded in time, but  travel
slower, may exhibit recurrent phenomena and do not look like scattering states  
(see, for example, \cite{Last} and discussion in 
\cite{CK3} for more detailed explanation and further references). A significant effort in quantum 
mechanical scattering theory was devoted to proving asymptotic completeness, and thus 
 absence of singular continuous spectrum, in a wide variety of situations (see, e.g. 
\cite{CFKS} for a review of some examples). However, to the best of our knowledge, there are no 
examples of Schr\"odinger operators for which wave operators exist but are not asymptotically 
complete due to the presence of singular continuous spectrum. In the problem that
we consider, a classical result is the existence of  wave operators for $V \in L_1,$ and in this class 
 wave operators are also asymptotically complete. However recently in \cite{CK3}, the
existence of (modified) wave operators for $V \in L_p,$ $p<2$ was established. Moreover, usual M\"oller wave operators
\eqref{wo} exist if $\lim_{x \rightarrow \infty} \int_0^x V(y)\,dy$ exists. 
Therefore, an example of (conditionally integrable) potential in $L^p,$ $p<2$ with imbedded singular   
continuous spectrum would also provide the first known example of a Schr\"odinger operator for 
which wave operators exist but fail to be asymptotically complete due to the presence of the singular 
continuous spectrum.

Our main goal in this paper is to prove the following 
\begin{theorem}\label{main}
Given any positive function $h(x)$ tending to infinity as  $x$ grows,  
there exist potentials $V(x)$ such that $|V(x)| \leq \frac{h(x)}{1+x}$ 
and the operator $H_V$ has non-empty singular continuous spectrum.
\end{theorem}
 It is easy to see that without loss of generality one can assume that $h(x)$
is monotone.  For the rest of the paper, we are going to operate under this assumption, 
which is convenient for technical reasons. 
We are going to provide a fairly explicit construction of such potentials. 
As a corollary of such construction we are also going to prove 
\begin{theorem}\label{asin}
Given any positive function $h(x)$ tending to infinity as  $x$ grows,  
there exist potentials $V(x)$ such that $|V(x)| \leq \frac{h(x)}{1+x},$ 
and the wave operators for the operator $H_V$ exist but are not asymptotically 
complete. 
\end{theorem} 
 
Our examples can be easily generalized to any dimension by looking at spherically 
symmetric potentials.

Let us outline briefly the idea of construction. Controlling singular continuous 
spectrum is hard because unlike point spectrum, one has to control in a certain sense
an uncountable set of energies, and there is no simple criteria (like $L^2$) 
for what  behavior of solutions implies existence of singular continuous component. 
It turns out that one 
can reduce matters to keeping track on each step of only finite number of energies, but the price 
one has to pay is having to control rather precisely norms of the solutions at these energies.
The control has to be much more precise than in constructions of Naboko and Simon. 
At the heart of the construction is a lemma which allows instead of one imbedded eigenvalue in 
Wigner-von Neumann-type construction to get two, each having an eigenfunction (normalized at the
origin) with the $L^2$ norm squared exactly twice the $L^2$ norm squared 
of the eigenfunction in a usual Wigner-von Neumann construction in a certain asymptotic regime. 
The $L^2$ norm squared of the eigenfunction
is the inverse of the point weight the measure $\mu$ assigns to a given eigenvalue. Thus the lemma
will allow to ``spread  out" evenly the weight of one point mass to two nearby points. 
This allows to carry out Cantor-like construction, building a sequence of potentials which have point 
measure components approximating singular continuous measure. In the limit, one gets the result. 

The methods of this paper can be generalized to the discrete Schr\"odinger operators, by using the discrete analog of 
the Pr\"ufer transform (see, e.g. \cite{KLS}). The results are parallel to the continuous case.
We plan to address these applications elsewhere.

\section{A Splitting Lemma}\label{split1}

For the rest of the paper, let us fix for simplicity some boundary condition 
at the origin, for example Dirichlet. Let us denote by $\mu$ the spectral measure corresponding  
to this boundary condition. 
Throughout the paper, we are going to use notation $C$ for universal
constants (not necessarily the same), more precisely, for the constants which do not depend in any way 
on the step of the inductive construction. 
A  starting point is the following elementary observation. 
\begin{lemma}\label{el}
Let $E$ be an eigenvalue of $H_V,$ and $u(x,E)$ be the corresponding eigenfunction
normalized by the condition $u'(0,E)=1.$ Then $\mu(E)= \|u(x,E)\|_{L^2(\reals^+)}^{-2}.$ 
\end{lemma}
\begin{proof}
Denote $v(x,E)$ the solution satisfying orthogonal boundary condition $v(x,E)=1,$ $v'(x,E)=0.$ 
Recall that the Weyl $m$-function is defined for $z \in \complex \setminus \reals$ as a unique complex number $m(z)$ 
such that $v(x,z)+m(z)u(x,z) \in L^2$ (assuming we are in the limit point case). 
We have  
(see e.g. \cite{CL}, Chapter 9, Section 3)
\begin{equation}\label{mrep}
 m(z)-m(z_0)  = \int\limits_\reals \left( \frac{1}{t-z} -\frac{1}{t-z_0} \right) d\mu(t),
\end{equation}
where $z_0$ is some fixed complex number.
From \eqref{mrep} it follows that 
\[ \mu(E)= \lim_{\epsilon \rightarrow 0} \epsilon |m(E+i\epsilon)| = 
\lim_{\epsilon \rightarrow 0} \epsilon \Im m(E+i\epsilon). \]
On the other hand, 
\[  \Im m(E+i\epsilon) = \epsilon \|v(x,E+i\epsilon)+m(E+i\epsilon) u(x,E+i\epsilon)\|^2_{L^2(\reals^+)} \]
(see, e.g. \cite{CL}).
Multiplying both sides of the above equality by $\epsilon$ and passing to the limit $\epsilon \rightarrow 0,$ 
we get 
\[ \mu(E) = \mu(E)^2 \|u(x,E)\|^2_{L^2},  \]
which is exactly the claim  of the lemma.
\end{proof}

Given the energy $E=k^2,$ it will be convenient to introduce Pr\"ufer variables 
$R$ and $\theta,$ $R^2 = (u')^2 + k^2 u^2$ and $\theta = \tan^{-1}(ku/u').$
 Then it is easy to see that 
\begin{eqnarray}\label{pruferam}
(\log R(x,k)^2)' = \frac{1}{k}V(x) \sin 2\theta(x,k) \\
\label{pruferan}
\theta(x,k)' = k - \frac{1}{k} V(x) (\sin \theta)^2. 
\end{eqnarray}

There are two essential new ideas which play important role in the result: the Splitting Lemma, and a
relation between $L^2$ norm of $R$ and derivative of the Pr\"ufer angle $\theta$ in energy.
The purpose of this section is to prove the first key lemma which is at the heart of construction,
the Splitting Lemma. To clarify its meaning, let us formulate first 
a simplified version of this result. Fix an arbitrary positive energy $k_0$ (we are going to work 
with values of quasimomentum $k$ throughout the paper, 
but still call those ``energies" for simplicity). 
Fix a parameter $g,$ which we think of as large. Define $\tilde{g}=\exp(g^{3/4}).$ Fix another 
parameter, $f,$ which we think of as small. (It is useful to keep in mind that in the actual 
construction, we will 
be able to choose $g$ as large as we want and then after that to choose $f$ as small as we want). 
Let $\delta k = \tilde{g}^{-1} f,$ 
and denote $k_1,$ $k_2,$ $k_2>k_1,$ the ends of the interval of length $\delta k$ centered at $k_0.$ 
Denote $\theta_i (x)= \theta(x,k_i)$ and $\delta \theta = \theta_2-\theta_1.$ 
Furthermore, define a small angle $\alpha$ by a condition $f \sin \alpha = \delta k.$ 
We are going to use notation 
\[ x_\beta = {\rm min} \{ x: \delta \theta (x) = \beta \}. \]
Define a potential $V(x)$
on $[0,\infty)$ by 
\begin{equation}\label{splipot}
V(x) = \left\{ \begin{array}{ll} -2 f k_0 \sin(\theta_1+\theta_2), & 0 \leq x \leq x_{\pi-\alpha} \\
2f k_0 \sin(\theta_1 +\theta_2), & x_{\pi-\alpha} \leq x \leq f^{-1}g \\
-{\rm min} \left( \delta k, \frac{g}{2x} \right) \frac{k_0}{4} (\sin 2\theta_1+\sin 2\theta_2), &
x > f^{-1}g. \end{array} \right.
\end{equation}
The definition of $V$ involves the Pr\"ufer angles $\theta_i,$ and $\theta_i$ in their turn depend on $V$ through 
\eqref{pruferan}. This ambiguity is easily resolved, however, by substituting the expression \eqref{splipot}
for $V$ into \eqref{pruferan} and solving the resulting nonlinear system for $\theta_i.$ It is not difficult
to see from \eqref{pruferan} that this resulting system has unique piecewise smooth global solution by standard ODE 
existence and uniqueness theorem. The functions $\theta_i(x)$ are smooth apart from three points where $V$ is not smooth,
and where $\theta_i'$ (or $\theta_i''$) may jump. 
We can then define $V$ in terms of these $\theta_i$ by \eqref{splipot}
recipe, and then solutions of \eqref{pruferan} will coincide with $\theta_i$ by uniqueness.
For the rest of the paper, we are going to consider on each step $V$ defined in terms of Pr\"ufer
angles at a finite number of energies without further explanation.  
We are going to ignore for now the issue whether $V$ is well-defined in a sense that $x_{\pi-\alpha}<f^{-1}g$;
this will be verified later for all $g$ large enough. We notice that the constant $f$ in \eqref{splipot} can be replaced 
by an appropriate slowly varying function with no difference for the result. We choose a particular representation \eqref{splipot}
to simplify technicalities.  
Then we have 
\begin{lemma}\label{sp1}
Assume that $g$ is sufficiently large, and $f$ is sufficiently small.
Then for a potential $V(x)$ given by \eqref{splipot}, both $k_1$ and $k_2$ are eigenvalues.
Moreover,
\begin{equation}\label{spnoes}
\|R(x,k_i)\|^2_{L^2(\reals^+)} = 2f^{-1}(1+O(f e^g, \tilde{g}^{-1/2} g))
\end{equation}
provided that $R(x,k_i)$ are normalized by a condition $R(0,k_i)=1,$ $i=1,2.$
\end{lemma}
We are not going to prove Lemma~\ref{sp1}, but directly a more advanced version which is 
going to be needed in the construction. The difference is that we will need to handle not just a pair
of eigenvalues, but $2^{n-1}$ pairs of eigenvalues on the $n$th step simultaneously. 
However, we are going to use this more transparent formulation to make a few clarifying remarks. 
First of all, observe that $|V(x)|  \leq h(x)/(x+1)$ provided that $f$ is sufficiently small. 
Indeed, we basically need to ensure that
\[ 2k_0 f \leq h(f^{-1}g)/(f^{-1}g+1), \]
or $2k_0 (g+f) \leq h(f^{-1}g).$ This clearly holds if for fixed $g$ we take $f$ sufficiently small.
Next, consider a Wigner-von Neumann-type potential $V_0(x)$ defined by 
\[ V_0(x) = \left\{ \begin{array}{ll} -2k_0 f \sin 2\theta_0, & 0 \leq x \leq f^{-1}g \\
  -{\rm min} \left( \delta k, \frac{g}{2x} \right) \frac{k_0}{2} \sin 2\theta_0, &
x > f^{-1}g. \end{array} \right.\]
It is fairly straightforward to check using \eqref{pruferam}, \eqref{pruferan}
that $V_0$ leads to an eigenvalue at $k_0$ 
with $\| R(x,k_0) \|^2_{L^2(\reals^+)} = f^{-1} (1+ O(f g, e^{-g})).$ 
Moreover, one cannot get a smaller (in the main term) norm for an eigenvalue at $k_0$ with a potential
$V$ satisfying $|V(x)| \leq f.$ By Lemma~\ref{el}, this means that 
we cannot have a point mass
bigger than $f(1+O(fg,e^{-g}))$ in the spectral measure for such potentials.
Thus Lemma~\ref{sp1} provides a way, using potential satisfying the same upper bound as $V_0$, 
to split the eigenvalue $k_0$ into two nearby eigenvalues
with equal weights without loss of the total weight (in the asymptotic regime). 
The fact that it can be done is not surprising, but the construction is not trivial.
For example, taking just $V=\frac{1}{2}(V_1+V_2),$ where $V_i$ is given by the same expression as $V_0$ 
at $k_i,$ seems to lead to a loss of a constant factor in the total mass of two eigenvalues and thus 
does not work for the Splitting Lemma. 
Moreover, one cannot hope for a coefficient smaller than $2$ in front of $f^{-1}$ in 
\eqref{spnoes}. As the construction below will show, that would have led by splitting 
to spectral measures giving arbitrary  large weight to a fixed finite interval. But the spectral measure
$\mu$ is well known to satisfy $\int \frac{d\mu(t)}{1+|t|^2} <C,$ with $C$ uniform for uniformly bounded
$V.$ Therefore, the factor $2$ in the Splitting Lemma is optimal. 
 

After these preliminary remarks, let us consider the general case. 
On the $n$th step of our construction, we are going to look at $2^n$ energies 
$k^{(n)}_j,$ $j=1, \dots 2^n,$ ordered in increasing order. 
We are going to assume that for every $n,$ all these energies 
lie in a fixed compact interval away from zero; moreover, we are going to assume 
that 
\begin{equation}\label{genas}
2{\rm max}_{j,j'}|k_j^{(n)}-k_{j'}^{(n)}| \leq {\rm min} (k_j^{(n)})
\end{equation} 
for every $j,$ $n.$ This assumption is made for technical convenience and will be easy to satisfy 
on each step of the inductive construction. 
We denote $\theta^{(n)}_j= \theta (x, k ^{(n)}_j),$
$R^{(n)}_j= R(x, k^{(n)}_j),$ 
$\delta \knj = k^{(n)}_{2j}-k^{(n)}_{2j-1},$ and $\delta \tnj = 
\theta^{(n)}_{2j}-\theta^{(n)}_{2j-1}.$ Then $V^{(n)}(x)$ is going to be defined 
by $V^{(n)}(x) = \sum_{j=1}^{2^{n-1}} V^{(n)}_j(x),$ where 
\begin{equation}\label{sppo}
V^{(n)}_j(x) = \left\{ \begin{array}{ll} -2 f_j^{(n)} k_j^{(n-1)} \sin(\theta^{(n)}_{2j-1}+\theta_{2j}^{(n)}),
 & 0 \leq x \leq x^{(n)}_{j,\pi-\alpha^{(n)}_j} \\
2 f_j^{(n)} k_j^{(n-1)} \sin(\theta^{(n)}_{2j-1}+\theta_{2j}^{(n)}), 
& x^{(n)}_{j,\pi-\alpha^{(n)}_j} \leq x \leq (f^{(n)}_j)^{-1}g \\
-{\rm min} \left( \delta \knj, \frac{g_n}{2x} \right) \frac{k_j^{(n-1)}}{4} (\sin 2\theta_{2j-1}^{(n)}
+\sin 2\theta_{2j}^{(n)}), &
x > (f_j^{(n)})^{-1}g. \end{array} \right.
\end{equation}
Here $\alpha_j^{(n)}$ is defined by $f_j^{(n)} \sin \alpha_j^{(n)} = \delta k_j^{(n)},$ and 
\[ x^{(n)}_{j,\beta}  = {\rm min} \{ x: \delta \theta^{(n)}_j (x) = \beta \}. \]
Notice that the value of $g_n$ is going to be the same for all $j.$ Also, as before, we will denote 
$\tilde{g}_n = \exp(g_n^{3/4}).$ The parameters $g_n,$ $f_j^{(n)}$ and $\delta k_j^{(n)}$ are to 
be chosen on the $n$th step of construction. The energies $k_j^{(n-1)}$ are given (from the previous step). 
On the $n$th step, each $k_j^{(n-1)}$ splits into two eigenvalues $k_{2j-1}^{(n)},k^{(n)}_{2j};$ that is, the interval 
   $(k_{2j-1}^{(n)},k^{(n)}_{2j})$ is centered at  $k_j^{(n-1)}$ and has length $\delta k_j^{(n)}.$
Throughout the paper, we will assume the relationship 
$ f_j^{(n)} = \tilde{g}_n \delta k_j^{(n)}. $
In the actual construction, $V^{(n)}(x)$ will be given by \eqref{sppo} only to the right of some value of $x=x_n,$ 
and $R(x_n,\knj) \ne 1.$ However in this section we will consider $V^{(n)}(x)$ defined by \eqref{sppo}
on $[0,\infty)$ and we will assume that $R(0,k) = 1$ and $\tnj (0,k) =0$ for all $k.$  
Let us introduce one more parameter that we need, $a_n = {\rm min}_{j,j'} |k^{(n-1)}_j-k^{(n-1)}_{j'}|$ 
(essentially, in the construction process we will have $a_n = {\rm min}_j (\delta k^{(n-1)}_j)$).
To keep notation compact, henceforth in this section  we are going to omit the index ``$n$" for most variables,
including $V$
(but with an exception of $k$'s, where members of two different steps participate explicitly in a construction). 
\begin{lemma}[Splitting Lemma]\label{split}
Let $V(x)$ be given by \eqref{sppo}, and assume that 
\begin{equation}\label{assplit}
4\sum_j f_j k^{(n-1)}_j < \frac13 a,\,\,\,\,\, g >>1,\,\,\,2^n g a^{-2} \sum\limits_{j=1}^{2^{n-1}} f_j<<1, 
 \,\,\,\,\, \delta k_j < \frac{1}{12} a, 
\end{equation} 
for $j=1, \dots, 2^{n-1}.$ 
Then for any $l=1, \dots , 2^{n-1},$ 
\begin{equation}\label{splitting}
\|R(x,k^{(n)}_{2(l-1)+i})\|^2_{L^2(\reals^+)} = 2f_l^{-1} \left(1+ O\left( 2^n e^{g}a^{-2} \sum\limits_{j=1}^{2^{n-1}} f_j, 
 \tilde{g}^{-1/2} g\right)\right), 
\end{equation}
where $i=1,2.$ 
\end{lemma}
\noindent \it Remark. \rm The condition $g>>1$ means that $g$ needs to be greater than some universal 
constant, the value of which can be derived from the proof. The same interpretation applies to the other condition in 
\eqref{assplit} involving ``$<<$''. 
 Notice also that the fourth condition
basically follows from the first two and $\delta k_j = \tilde{g}^{-1} f_j$, 
but we state it separately for convenience. \\

We are going to fix $l=1,$ $i=1$ and prove the result for this case; other $l,i$ are treated in exactly the same way. 
Write $V=V_1+W,$ where $V_1= V^{(n)}_1$ and $W$ is the rest of the potential. Direct substitution of \eqref{sppo} into
\eqref{pruferam}, \eqref{pruferan} and some trigonometry leads to the following equations valid for 
$x < f_1^{-1}g:$ 
\begin{equation}
\label{keyR} (\log R^2(x,k_1^{(n)})'  =  -\frac{k^{(n-1)}_1}{k^{(n)}_1}\tilde{f}_1(x) 
(\cos \delta \theta_1 - \cos(\theta_2+3\theta_1))
+ \frac{W(x)}{k^{(n)}_1}\sin 2\theta_1
\end{equation}
\begin{eqnarray}
\label{keyan}
(\delta \theta_1)' & = & \delta k_1 \left( 1- \frac{V(x)}{2k^{(n)}_1 k_2^{(n)}}(1+\cos 2\theta_2) \right) + \\
\nonumber &  & + \left( \frac{k^{(n-1)}_1}{k^{(n)}_1}\tilde{f}_1(x)  (1-\cos2(\theta_1+\theta_2)) -
\frac{W(x)}{k_1^{(n)}} \sin(\theta_1+\theta_2) \right) \sin \delta \theta_1; 
\end{eqnarray}
here $\tilde{f}_1(x)=\pm f_1$ with a change from $+$ to $-$ at $x = x_{1,\pi - \alpha_1}.$  
For the rest of this section, we are also going to omit the index ``$1$'' in $x_{1,\beta}$ and write 
simply $x_\beta$ for $x$ such that $\delta \theta_1 (x_\beta) =\beta.$ \\

\noindent \it Remark. \rm Equations \eqref{keyR} and \eqref{keyan} are fairly complicated, but not all terms are of equal
importance. At the heart of the matter is the following simpler nonlinear dynamical system:
\begin{eqnarray}\label{easysys}
(\log R^2)' = -\tilde{f}_1(x) \cos \delta \theta_1 \\
(\delta \theta_1)' = \delta k_1 + \tilde{f}_1(x) \sin \delta \theta_1. \nonumber
\end{eqnarray}
In some sense, we will show that the rest of  
terms in \eqref{keyR}, \eqref{keyan} produce only small corrections to the behavior of this system. 
It is quite instructive (and technically simple) to prove a version of the Splitting Lemma for the system
\eqref{easysys}, but to save the space, we are going to treat directly the system \eqref{keyR}, \eqref{keyan}. \\ 

From \eqref{keyan} and $\delta \theta_1(0)=0,$ we can write
\begin{equation}\label{andif}
\delta \theta_1 =  \delta k_1 \int\limits_0^x  
e^{\int_y^x \xi(s) \,ds} (1+ O(|V(y)|))\,dy,
\end{equation}
where 
\begin{equation}\label{xif}
 \xi (y) = \frac{k^{(n-1)}_1}{k^{(n)}_1}\left( \tilde{f}_1(y) 
(1-\cos 2(\theta_1 + \theta_2)) - \frac{W(y)}{k_1^{(n-1)}}\sin(\theta_1+\theta_2)
\right) \frac{\sin \delta \theta_1}{\delta \theta_1}. 
\end{equation}

We need the following technical lemma 
\begin{lemma}\label{aux1}
Under the assumptions of the Splitting Lemma, we have for $y<x \leq f_1^{-1}g$ 
\begin{equation}\label{poest}
 \int\limits_y^x \xi (s) \,ds = \int\limits_y^x \tilde{f}_1(s) \frac{\sin \delta \theta_1}{\delta \theta_1}\,ds
(1+O(\delta k_1)) +  O(\scripte_n), 
\end{equation} 
where $\scripte_n=2^n (\sum_j f_j) g a^{-2}.$
\end{lemma}
\begin{proof}
Let us consider one of the  terms entering into $\xi,$ 
\begin{equation}\label{ont}
 \int\limits_y^x V_j(s)  \sin(\theta_1+\theta_2) \frac{\sin \delta \theta_1}{\delta \theta_1}\,ds, 
\end{equation}
where $2 \leq j \leq 2^{n-1}$ is arbitrary. 
Recall that $V_j(s) = \tilde{f}_j(s) \sin(\theta_{2j-1}+\theta_{2j})$ for $x \leq (f_j)^{-1}g,$ where
$\tilde{f}_j = \pm f_j$ with only at most one jump (we do not assume apriori that 
$x_{j,\pi-\alpha_j} < (f_j)^{-1}g,$ but will show this below). For $x > (f_j)^{-1}g,$ 
$V_j(s) = \tilde{f}_j(s) (\sin(2\theta_{2j-1} )+\sin (2\theta_{2j})),$ where $|\tilde{f}_j(s)| \leq (k_j^{(n-1)}g)/2s,$ 
$|\tilde{f}_j'(s)| \leq (k_j^{(n-1)}g)/(2s^2).$ Notice that by assumption \eqref{assplit} we have 
\begin{equation}\label{ansep3}
 |(\theta_{2j-1}+\theta_{2j}-\theta_1-\theta_2)'| \geq \frac{a}{2} , \,\,\,\,
 |(2\theta_{2j-1}-\theta_1-\theta_2)'| \geq \frac{a}{2}, \,\,\,\, |(2\theta_{2j}-\theta_1-\theta_2)'| \geq \frac{a}{2}. 
\end{equation}
Substituting the expression for $V_j$ into \eqref{ont} and using formula for the product of sines, we reduce 
matters to estimation of the integrals of type
\[ I_\pm(x,y)=\int\limits_y^x \tilde{f}_j  \cos(2\theta_{2j}\pm (\theta_1+\theta_2)) 
\frac{\sin \delta \theta_1}{\delta \theta_1}\,ds.\] 
Consider the ``$-$" case, and integrate by parts integrating 
\[ \cos(2\theta_{2j}-(\theta_1+\theta_2))(2\theta_{2j}-(\theta_1+\theta_2))'. \]
Observing that $|\theta_j''(s)| \leq C (|V(s)|+|V'(s)|)$ for any $j,$ 
$|(\delta \theta_1)'(s)|  \leq C(|V(s)|+\delta k_1)$ 
and using the inequality \eqref{ansep3} for the derivative of the argument of cosine, we obtain 
\[ |I_-(x,y)| \leq C a^{-1} \left( f_j  + \int\limits_y^x  \left( |\tilde{f}_j'(s)|  + |\tilde{f}_j(s)| (|V(s)|+\delta k_1) +
|\tilde{f}_j(s)| (|V(s)|+|V'(s)|) a^{-1} \right)ds \right), \]
  where $C$ is a universal constant. The discontinuities of $\tilde{f}_j$ are taken into account in off-integral term.
Since $|V(s)|,|V'(s)| \leq C \sum_l |\tilde{f}_l(s)|,$ $\delta k_1 < f_1$ and $a \leq C,$ 
we arrive at a bound valid for all $y<x \leq f^{-1}_1 g$
\[ |I_-(x,y)| \leq C\left( f_j a^{-1} +  a^{-2} \int\limits_y^x |\tilde{f}_j(s)|\sum_l |\tilde{f}_l(s)|\,ds \right). \]
The estimate on $I_+$ is similar, but involves a constant independent of $n$ instead of $a$ (arising from the 
minimal possible size of $k^{(n)}_j$). 
The contribution of the term $\tilde{f}_1(s) \cos(2(\theta_1+\theta_2))$ in \eqref{xif} is bounded similarly 
to $I_+.$ Summing up all bounds, we obtain at most  
\begin{equation}\label{err1}
 C \left( \sum_j f_j  a^{-1} + a^{-2} \int\limits_y^x \left(\sum_j |\tilde{f}_j(s)|\right)^2 \,ds \right). 
\end{equation}
 We can compute explicitly that $\|\tilde{f}_j\|_{L^2}^2 \leq 2f_j g,$ and therefore the total bound does not exceed
$ C a^{-2}2^n (\sum_j f_j)g.$ Notice that the first term in \eqref{err1} is included into this error since $a$ is less 
than a fixed constant and $g>1.$ 
This is exactly the error term in \eqref{poest}; the origin of the main term is clear.  
\end{proof}
\noindent \it Remark. \rm The error $O(\delta k_1)$ in the main term of \eqref{poest} gives after integration 
$O(g\delta k_1)=O(\tilde{g}^{-1}gf_1)$ which can also be subsumed into $\scripte_n.$   

Next we show that, provided
that the error terms in \eqref{poest} are small, 
we have $x_{\pi-\alpha_1} < f_1^{-1}g,$ and thus our potential is well defined.  
\begin{lemma}\label{el2}
For sufficiently  large $g$ and for sufficiently small $\scripte_n,$ we have 
\begin{equation}\label{xpial}
x_{\pi-\alpha_1} \leq \pi f_1^{-1} \log (2 \tilde{g})(1+C\scripte_n).
\end{equation}   
\end{lemma}
\noindent \it Remark. \rm Similarly to Lemma~\ref{split}, by 
sufficiently large (or small) we mean that there exist universal 
constants such that if the value of $g$ (respectively $\scripte_n$) exceeds (is less than) these constants, 
the result holds. 
 
\begin{proof}
Let us first estimate from above $x_{\pi/2}.$ On the interval $[0, x_{\pi/2}]$ we have 
\[ 2/\pi \leq (\sin \delta \theta_1)/\delta \theta_1 \leq 1. \]
Thus, according to \eqref{andif},  Lemma~\ref{aux1} and the remark after its proof, we have 
\begin{equation}\label{thineqlow}
 \delta \theta_1 (x) \geq e^{\frac{2}{\pi} f_1 x}  \frac{\pi \delta k_1}{2 f_1}( 1 - e^{-\frac{2}{\pi}f_1x} )
(1 -C\scripte_n)
\end{equation}
for $x < (f_1)^{-1}g$ and such that $\delta \theta_1 \in [0, \pi/2].$
It is clear from \eqref{thineqlow} that if $g$ is large enough, then $x_{\pi/2} <f_1^{-1}g.$ 
Indeed, recall that $\delta k_1/f_1 = e^{-g^{3/4}}.$ If $f_1^{-1}g$ were less than $x_{\pi/2},$ 
 \eqref{thineqlow} would imply that 
\[ \delta \theta_1 ( f_1^{-1}g) \geq e^{\frac{2}{\pi} g -g^{3/4}} \frac{\pi}{2} \left( 1- e^{-\frac{2}{\pi}g}\right)
(1-C\scripte_n), \]
an obvious contradiction for $g$ large enough and sufficiently small $\scripte_n.$
Moreover,  provided that $\scripte_n$ is small, 
we can estimate (assuming also $x_{\pi/2} \geq \frac{\pi}{2}(f_1)^{-1} \log 2,$
or else we have a nice upper bound)
\[ e^{\frac{2}{\pi} f_1 x_{\pi/2} } \leq 2 \tilde{g}\left(1+ C \scripte_n \right). \]
From this inequality it is easy to see that
\begin{equation}\label{xpi2} 
 x_{\pi/2} \leq \frac{\pi}{2} f_1^{-1} \log(2 \tilde{g})(1+C\scripte_n).
\end{equation}
Imagine for a moment that $\tilde{f}_1(x) = f_1$ for $x \leq x_{\pi}.$ Notice that under a change of 
variables $\overline{\delta \theta}_1= \pi - \delta \theta_1,$ $\overline{x} = x_{\pi} -x$ the equation 
\eqref{keyan} transforms into the equation for $\overline{\delta \theta}_1$ of the same form as \eqref{keyan},
with the same initial condition and same estimates available. In particular, the point $\overline{x}_{\pi/2}$
such that $\overline{\delta \theta}_1(\overline{x}_{\pi/2}) = \pi/2$ satisfies the same bound as $x_{\pi/2}.$ 
But $\overline{\delta \theta}_1 (\overline{x}) = \pi - \delta \theta_1 (x_{\pi} -x),$ and thus 
$\overline{x}_{\pi/2} = x_\pi -x_{\pi/2}.$ 
But $x_{\pi-\alpha_1} \leq x_\pi,$ and the above observation together with \eqref{xpi2} finishes the proof 
of the lemma.
\end{proof}

Next, we prove the following essential
\begin{lemma}[Small Angle Lemma]\label{sman} 
Let $\gamma$ be a small angle, in particular such that $\gamma < g^{-1}.$ 
Denote $x_\gamma$ the smallest value of $x$ where $\delta \theta_1(x) = \gamma.$ 
Then 
\begin{equation}\label{declev}
R(x_\gamma, k_1^{(n)})^2 =  \frac{\delta k_1}{\gamma f_1+\delta k_1} \left(1+O(\gamma^2 g, 
\scripte_n )\right).
\end{equation} 
Moreover,
\begin{equation}\label{normcon}
\|R(x,k^{(n)}_1) \|^2_{L^2(0,x_\gamma)} = f_1^{-1} \frac{1}{1+\tilde{g}^{-1}\gamma^{-1}}\left(1+O(\gamma^2 g,  
\scripte_n)\right). 
\end{equation}
\end{lemma}
\begin{proof}
Notice that on $[0,x_\gamma],$ we have $(\sin \delta \theta_1)/\delta\theta_1 = 1+O(\gamma^2).$ 
Obesrve that by Lemma~\ref{el2} we have that $x_\gamma < f_1^{-1}g.$ 
By \eqref{andif}, Lemma~\ref{aux1} and the definition of potential we have that 
\[ \delta k_1 \int\limits_0^{x_\gamma}  e^{f_1 (x_\gamma - y) (1+O(\gamma^2)) +O(\scripte_n)}
(1+O(|V(y)|)\,dy  = \gamma. \]
Thus, we obtain 
\begin{equation}\label{deccon}
e^{f_1 x_\gamma} (1- e^{-f_1 x_\gamma}) f_1^{-1} \delta k_1 (1+O(\scripte_n, \gamma^2 g)) =\gamma.
\end{equation}
We included the $O(|V(y)|)$ error into $\scripte_n,$ since we can assume freely 
that $a<1$ and $g>1.$ Thus 
\begin{equation}\label{addb34}
e^{f_1 x_\gamma} = \left( \frac{\gamma f_1}{\delta k_1} +1 \right) \left( 1+ O(\gamma^2 g, \scripte_n)\right). 
\end{equation} 
Now in the equation \eqref{keyR} for the amplitude, we have for $x < f_1^{-1}g$
\begin{equation}\label{keyRest}
\log( R(x, k_1^{(n)})^2) = -\int\limits_0^x \tilde{f}_1(y) \cos \delta \theta_1(y)\,dy + O(\scripte_n),
\end{equation}
by an estimate directly analogous to the estimates of Lemma~\ref{aux1} (but easier).
Therefore, for $x \leq x_\gamma,$ we have 
\begin{equation}\label{RR1}
 R(x, k_1^{(n)})^2 = e^{-f_1 x}(1+ O(\gamma^2 g, \scripte_n)). 
\end{equation}
The first statement of the Lemma, estimate \eqref{declev}, follows directly from \eqref{RR1} and \eqref{addb34}.
Also, the $L^2$ norm of $R(x,k_1^{(n)})$ on $[0,x_\gamma]$ satisfies
\[ \|R(x, k_1^{(n)})\|^2_{L^2(0,x_\gamma)} = f_1^{-1} (1- e^{-f_1 x_\gamma})(1+O(\gamma^2 g, \scripte_n))=
f_1^{-1}\frac{1}{1+\tilde{g}^{-1} \gamma^{-1}}(1+O(\gamma^2 g,\scripte_n)). \]
In the last step we used \eqref{declev} and the fact that $\delta k_1 = \tilde{g}^{-1}f_1.$
\end{proof}

Next we need a lemma describing the solution of \eqref{keyR}, \eqref{keyan} on the interval $[0,x_{\pi-\alpha_1}].$ 
Let $R(x,k^{(n)}_1),$ $\delta \theta_1$ solve \eqref{keyR}, \eqref{keyan} and assume that $\tilde{R}(x),$ $\tilde{\delta \theta}_1(x)$ 
solve a variant of a system \eqref{easysys}
\begin{eqnarray}\label{simpsysR}
(\log \tilde{R}^2(x))' & = & -f_1 \cos \tilde{\delta \theta}_1(x) \\
(\tilde{\delta \theta}_1(x))'& =& \delta k_1 + f_1 \sin \tilde{\delta \theta}_1(x) \label{simpsysan}
\end{eqnarray}
for all $x.$
Assume furthermore that the initial conditions for $\tilde{R}$, $\tilde{\delta \theta}_1$ at $x=0$ are the same as for 
$R,$ $\delta \theta.$ Then the following lemma holds.
\begin{lemma}\label{idap}
For every $x \leq x_{\pi-\alpha_1},$ we have 
\begin{eqnarray}\label{tineq33}
|\tilde{\delta \theta}_1(x) - \delta \theta_1 (x)| \leq C \tilde{g}^{-1} \tilde{\scripte}_n \\
\label{Rineq33}
\tilde{R}^2(x) / R^2(x,k^{(n)}_1) = 1+ O(\tilde{\scripte}_n), 
\end{eqnarray}
where $\tilde{\scripte}_n = 2^n e^g a^{-2} \sum_j f_j.$ 
\end{lemma}
\begin{proof}
Denote $\delta \theta^*(x) = \tilde{\delta \theta}_1(x) - \delta \theta_1(x).$ 
From \eqref{keyan}, \eqref{simpsysan} it follows that for $x \leq  f_1^{-1}g,$
\begin{eqnarray*}
(\delta \theta^*)'(x) & = &\delta k_1 O(|V(x)|) + f_1 \frac{\sin \tilde{\delta \theta}_1 - \sin \delta \theta_1}{\delta \theta^*} 
\delta \theta^* + \\ && \left( f_1 \cos 2(\theta_1+\theta_2) + \frac{W(x)}{k_1^{(n)}} \sin (\theta_1+\theta_2) \right) \sin \delta \theta_1.
\end{eqnarray*}
Therefore, 
\begin{equation}\label{solan33}
\delta \theta^* (x) = \int\limits_0^x e^{\int_y^x \eta(s)\,ds} \left( \delta k_1 O(|V(y)|)+ 
\left( f_1 \cos 2(\theta_1+\theta_2) + \frac{W(x)}{k_1^{(n)}} \sin (\theta_1+\theta_2) \right) \sin \delta \theta_1 \right),
\end{equation}
where 
\[ \eta(s) = 2f_1 \frac{\sin \delta \theta^*}{\delta \theta^*} \cos (\delta \theta_1+\tilde{\delta \theta}_1). \]
Now all oscillatory terms in \eqref{solan33} (that is, all terms except the $\delta k_1 O(|V|)$ term) 
are estimated by integration by parts similarly 
to Lemma~\ref{aux1}. For example, for any of the $2^{n-1}-1$ terms entering into $W(y),$ we have for $x \leq x_{\pi-\alpha_1}$
\begin{eqnarray}\label{techest11}
 \left| \,\int\limits_{0}^x 
 e^{\int_y^x \eta(s)\,ds}  V_j(y) \sin(\theta_1 +\theta_2) \,dy \right| \leq \\
C a^{-1} \tilde{g} \left( f_j + \int\limits_0^x \left( |\tilde{f}_j'(y)|  + |\eta(y)| |\tilde{f}_j(y)| 
+ |\tilde{f}_j(y)| O(|V(y)|+|V'(y)|) a^{-1} \right) \,dy \right),  \nonumber
\end{eqnarray}
where 
\begin{equation}\label{auxf}
  \tilde{f}_j (y) = \left\{ \begin{array}{ll} -f_j, & 0 \leq x \leq x_{j,\pi-\alpha_j} \\
f_j, & x_{j,\pi-\alpha_j} \leq x \leq (f_j)^{-1}g \\
-{\rm min}(\delta k_j, \frac{g}{2x} ) \frac{k_j^{(n-1)}}{4}, & x > (f_j)^{-1}g. 
\end{array}
\right. 
\end{equation}
Lemma~\ref{el2} was used to estimate the exponential by $\tilde{g}.$ 
We make a  convention that the derivative $\tilde{f}_j'(y)$ does not involve $\delta$ functions from the two jumps 
of $\tilde{f}_j;$ those jumps contribute to the first, off-integral term in the estimate.
Notice that $|\eta(x)|,|V(x)|,|V'(x)| \leq C \sum_j \tilde{f}_j(x),$ where $C$ is a universal constant.
Thus the total bound in \eqref{techest11} does not exceed
\[ Ca^{-1} \tilde{g} \left(f_j + a^{-1} \int\limits_{0}^\infty \tilde{f}_j(y) \sum_l \tilde{f}_l(y) \,dy \right). \]
The term $f_1 \cos 2(\theta_1+\theta_2)$ is estimated similary (but gives smaller error), while 
the term $\delta k_1 O(|V(y)|)$ integrated over $[0, x_{\pi-\alpha_1}]$ with the exponential 
gives at most $(f_1)^{-1} \delta k_1 \tilde{g} \sum_j f_j$ by Lemma~\ref{el2}. Summing together all errors
we obtain 
\[ |\delta \theta^*| \leq C \tilde{g} \left( a^{-1} \sum_l f_l + a^{-2}
\int\limits_0^\infty \left( \sum_j \tilde{f}_j
\right)^2 \,dy  \right), \]
which leads to 
\begin{equation}\label{thofffin11}
|\delta \theta^*(x)| \leq C2^n  \tilde{g} g  a^{-2} \sum_j f_j
\end{equation} 
for $x \leq x_{\pi-\alpha_1}.$
This proves \eqref{tineq33}, since $e^g \geq Cg (\tilde{g}^2)$ for some universal $C$ and all $g.$  


By \eqref{keyR}, \eqref{simpsysR} and an argument parallel to Lemma~\ref{aux1},
\[ \log R^2(x, k_1^{(n)}) - \log \tilde{R}^2(x) = - \int\limits_0^x f_1 (\cos \delta \theta_1 - 
\cos \tilde{\delta \theta}_1 ) \,dy + O(\scripte_n). \]
Using \eqref{tineq33} and Lemma~\ref{el2} we obtain for $x \leq x_{\pi-\alpha_1}$ 
\begin{equation}\label{Rcomp33}
R^2(x, k_1^{(n)})/ \tilde{R}^2(x) = 1 + O( \tilde{\scripte_n}), 
\end{equation}
which finishes the proof (we used that evidently $\scripte_n < \tilde{\scripte_n}$ and so the former error is
absorbed by the latter). Also on the last step we generously exchanged $\tilde{g}^{-1}$ for $g^{3/4}$ to arrive at an error 
in \eqref{Rcomp33}.
\end{proof}

Now we are ready to proceed with the proof of Lemma~\ref{split}.
\begin{proof}[Proof of the Splitting Lemma]
The structure of the eigenfunctions obtained in Splitting Lemma is illustrated by Figure~\ref{struct}.

1. Consider first the interval $[0, x_{\pi/2}].$ The estimate for the norm of $R(x,k_1^{(n)})$ on the interval 
$[0,x_\gamma]$ for some small $\gamma$ is provided by Lemma~\ref{sman}.  By \eqref{declev} and \eqref{keyRest},
the norm on the interval $[x_\gamma, x_{\pi/2}]$  does not exceed 
$ (\frac{\gamma f_1}{\delta k_1}+1)^{-1} x_{\pi/2}(1+O(\gamma^2  g, \scripte_n)).$ 
Choosing $\gamma = \tilde{g}^{-1/2},$ for example, and using the bound \eqref{xpi2} on $x_{\pi/2},$ 
we find that 
\[ \| R(x, k_1^{(n)})\|_{L^2[x_\gamma,x_{\pi/2}]}^2 = f_1^{-1} O(\tilde{g}^{-1/2} g^{3/4}, \scripte_n). \]
Thus, altogether, 
\begin{equation}\label{xpi2norm}
\| R(x, k_1^{(n)})\|_{L^2[0,x_{\pi/2}]}^2 = f_1^{-1} (1+ O(\scripte_n,  \tilde{g}^{-1/2}g)). 
\end{equation}  

2. Next, we estimate norm on $[x_{\pi/2}, x_{\pi-\alpha_1}].$ 
Consider an auxiliary simplified system \eqref{simpsysR}, \eqref{simpsysan} of Lemma~\ref{idap}.
By \eqref{tineq33}, we have 
\[ \tilde{\delta \theta}_1( x_{\pi/2}) = \pi/2 +O(\tilde{g}^{-1} \tilde{\scripte_n}), \,\,\,
\tilde{\delta \theta}_1 (x_{\pi-\alpha_1}) = \pi -\alpha_1 + O(\tilde{g}^{-1} \tilde{\scripte_n}). \]
Denote $\tilde{x}_\beta = {\rm min} \{ x| \tilde{\delta \theta}_1 (x) = \beta \}.$ 
Notice that the system \eqref{simpsysR}, \eqref{simpsysan} is symmetric in a sense that for $x \leq \tilde{x}_\pi,$ 
\[ \tilde{\delta \theta}_1 (x) = \pi - \tilde{\delta \theta}_1 (\tilde{x}_\pi -x)\,\,\,{\rm  and}\,\,\, 
\tilde{R}^2(x) = \tilde{R}^2(\tilde{x}_\pi -x). \]
Observe also that $\alpha_1 = \tilde{g}_1^{-1} (1+ O(\tilde{g}^{-2})).$ 
Thus the norm 
\begin{equation}\label{noi33} 
\| \tilde{R} \|^2_{L^2[x_{\pi/2}, x_{\pi - \alpha_1}]} = 
  \| \tilde{R} \|^2_{L^2[\tilde{x}_{\alpha_1(1+O(\tilde{\scripte}_n))},\tilde{x}_{\pi/2(1+O(\tilde{\scripte}_n))}]}. 
\end{equation}
By identical arguments, Lemma~\ref{el2} and Lemma~\ref{sman} (Small Angle Lemma) apply to the system \eqref{simpsysR}, 
\eqref{simpsysan} (with the only difference that the conclusion does not involve the error $\scripte_n$ coming from 
oscillatory terms which are absent in the simplified system). 
Therefore, the considerations close to the first part of the proof above apply to the norm of $\tilde{R}$ and 
give
\begin{equation}\label{btR1}
 \|\tilde{R}\|^2_{L^2[0,\tilde{x}_{\pi/2(1+O(\tilde{\scripte}_n))}]}= f_1^{-1}
 (1+ O(\tilde{g}^{-1/2}g)). 
\end{equation}
By \eqref{normcon} of Lemma~\ref{sman}, we also have 
\begin{equation}\label{btR2}
 \|\tilde{R}\|^2_{L^2[0,\tilde{x}_{\alpha_1(1+O(\tilde{\scripte}_n))}]}= \frac{1}{2}f_1^{-1}
 (1+ O(\tilde{\scripte}_n, \tilde{g}^{-2}g)). 
\end{equation}
Therefore, the norm in \eqref{noi33} is equal to $\frac12 f_1^{-1}(1+ O(\tilde{\scripte}_n,\tilde{g}^{-1/2}g)).$
Now by the estimate \eqref{Rineq33} of Lemma~\ref{idap} we have 
\begin{equation}\label{xpi2al}
\| R(x, k_1^{(n)})\|_{L^2[x_{\pi/2}, x_{\pi-\alpha_1}]}^2 = \frac{1}{2}f_1^{-1} (1+ O(\tilde{\scripte}_n,  \tilde{g}^{-1/2}g)). 
\end{equation}
In addition, by Lemma~\ref{sman} and symmetry
\[ \tilde{R}^2(x_{\pi-\alpha_1}) = \tilde{R}^2(x_{\alpha_1(1+O(\tilde{\scripte_n}))})=
\frac{1}{2}(1+O(\tilde{\scripte}_n, \tilde{g}^{-2}g)). \]
By Lemma~\ref{idap}, this implies 
\begin{equation}\label{decal}
R(x_{\pi-\alpha_1},k_1^{(n)})^2 
= \frac{1}{2}(1+O(\tilde{\scripte}_n, \tilde{g}^{-2}g)). 
\end{equation}


\begin{figure}
\begin{center}
\scalebox{0.75}{\includegraphics{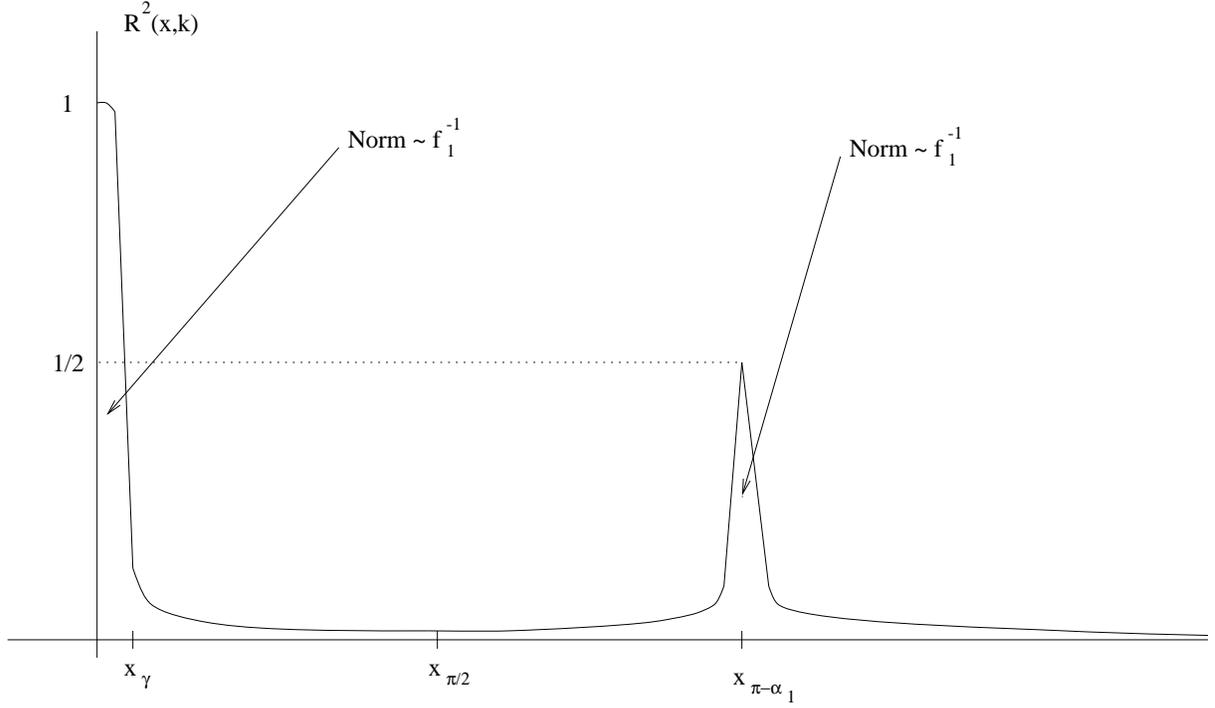}}
\caption{The structure of the eigenfunction at $k=k^{(n)}_{1,2}$. The error terms
are not pictured and lead to small oscillations around the main value.}
\label{struct}
\end{center}
\end{figure}

3. Next step is the consideration of $[x_{\pi-\alpha_1}, f_1^{-1}g]$ interval.
On this interval, in the idealized system \eqref{easysys} the angle $\delta \theta_1$ 
stays constant, and since it is very close to 
$\pi,$ this allows for a consistent decay in $R$ equation. The main goal therefore is to control the 
additional perturbation terms, which
are either small or oscillating. 
Define $\delta \overline{\theta}_1 = (\pi -\alpha_1) - \delta \theta_1,$ then
\[ (\delta \overline{\theta}_1)' = -\delta k_1 (1+O(|V(x)|)) + \sin(\delta \overline{\theta}_1 + \alpha_1)
\left( f_1(1-\cos 2(\theta_1+\theta_2)) + \frac{W(x)}{k_1^{(n)}} \sin (\theta_1+\theta_2)\right). \]
Using the relationship $\delta k_1 = f_1 \sin \alpha_1,$ we arrive at 
\begin{eqnarray} \nonumber
(\delta \overline{\theta}_1)' & = &  \left( f_1 (1-\cos 2(\theta_1 + \theta_2)) +
\frac{W(x)}{k_1^{(n)}} \sin (\theta_1 +\theta_2) \right) \times \\
\nonumber & & \times \left( \cos \alpha_1 \frac{\sin \delta \overline{\theta}_1}
{\delta \overline{\theta}_1}+ \sin \alpha_1 \frac{\cos \delta \overline{\theta}_1 -1}{\delta \overline{\theta}_1} \right)
\delta \overline{\theta}_1  - \\ & & \label{thetaoff} 
 \delta k_1 O(|V(x)|) + f_1 \sin \alpha_1 \cos 2(\theta_1+\theta_2)  - \frac{W(x)}{k_1^{(n)}} \sin \alpha_1
\sin(\theta_1+\theta_2),
\end{eqnarray}
with the initial data $\delta\overline{\theta}_1(x_{\pi-\alpha_1})=0.$ 
Denote $\zeta(x)$ the expression in front of $\delta \overline{\theta}_1$ in \eqref{thetaoff}. 
The solution to \eqref{thetaoff} 
can be written in the following form:
\begin{equation}\label{thetasol}
\delta \overline{\theta}_1(x) = \delta k_1 
\int\limits_{x_{\pi-\alpha_1}}^x 
 e^{\int_y^x \zeta(s)\,ds} \left( O(|V(y)|) -\cos 2(\theta_1+\theta_2) - \frac{W(y)}{f_1 k_1^{(n)}}
\sin(\theta_1+\theta_2) \right) \,dy. 
\end{equation}
Notice that for $x_{\pi -\alpha_1}  \leq y \leq x\leq (f_1)^{-1}g,$ 
\begin{equation}\label{zeta}
\int\limits_{y}^x \zeta(s)\,ds = f_1 \int\limits_{y}^x \left(\frac{\sin \delta \overline{\theta}_1}{\delta
\overline{\theta}_1}  \cos \alpha_1 + \frac{\cos \delta \overline{\theta}_1-1}{\delta \overline{\theta}_1} \sin \alpha_1 \right) \,ds +
O(\scripte_n)
\end{equation}
by the same argument as in the proof of Lemma~\ref{aux1}. 
Therefore, since $\alpha_1 \sim \tilde{g}^{-1},$ 
\begin{equation}\label{grow}
e^{ \int_{y}^x \zeta(s)\,ds} \leq e^{f_1 (x-y)}(1+O(\tilde{g}^{-1}g, \scripte_n))
\end{equation}
for $x \leq (f_1)^{-1}g.$ 
All oscillatory terms in \eqref{thetasol} (that is, all terms except the $O(|V|)$ term) are estimated by integration by parts similarly 
to Lemma~\ref{aux1} and Lemma~\ref{idap}. 
For example, for any of the $2^{n-1}-1$ terms entering into $W(y),$ we have for $x_{\pi-\alpha_1} \leq x \leq f_1^{-1}g$
\begin{eqnarray}\label{techest}
 \left| \,\int\limits_{x_{\pi-\alpha_1}}^x 
 e^{\int_y^x \zeta(s)\,ds}  V_j(y) \sin(\theta_1 +\theta_2) \,dy \right| \leq \\
C a^{-1} e^g \left( f_j + \int\limits_{x_{\pi-\alpha_1}}^\infty \left( |\tilde{f}_j'(y)|  + |\zeta(y)| |\tilde{f}_j(y)| 
+ |\tilde{f}_j(y)| O(|V(y)|+|V'(y)|) a^{-1} \right) \,dy \right),  \nonumber
\end{eqnarray}
where $\tilde{f}_j$ is given by \eqref{auxf} as before.
Similarly to the computation in Lemma~\ref{idap}, the total bound in \eqref{techest} is 
\[ Ce^g \left(f_j a^{-1} + a^{-2} \int\limits_{x_{\pi-\alpha_1}}^\infty \tilde{f}_j(y) \sum_l \tilde{f}_l(y) \,dy \right). \]
The term $\delta k_1 \cos 2(\theta_1+\theta_2)$ is estimated similary (but gives smaller error), while 
the term $\delta k_1 O(|V(y)|)$ integrated over $[x_{\pi-\alpha_1}, f_1^{-1}g]$ with the exponential 
gives at most $(f_1)^{-1} \delta k_1 e^g \sum_j f_j.$ Summing together all errors
we obtain 
\begin{equation}\label{thofffin}
|\delta \overline{\theta}(x)| \leq C2^n  \tilde{g}^{-1} g e^g a^{-2} \sum_j f_j. 
\end{equation} 
Now notice that 
\[ ( \log R(x,k_1^{(n)})^2 )' = \frac{k_1^{(n-1)}}{k_1^{(n)}} f_1 (-\cos(\delta \overline{\theta}_1 +\alpha_1)-
\cos(\theta_2+3\theta_1)) + \frac{W(x)}{k_1^{(n)}} \sin 2\theta(x, k_1^{(n)}). \]
Using \eqref{decal} and estimates parallel to the ones considered in Lemma~\ref{aux1} we find 
for $x_{\pi-\alpha_1} \leq x \leq (f_1)^{-1}g$ that  
\[ R^2(x,k_1^{(n)}) = \frac{1}{2} e^{-f_1(x-x_{\pi-\alpha_1})} \left(1+O(\tilde{\scripte}_n, g \tilde{g}^{-1}  )\right), \]
where as before
$ \tilde{\scripte}_n =2^n  e^g a^{-2} \sum_j f_j$
(we are being generous here, replacing $g^2 \tilde{g}^{-1} \tilde{\scripte}_n+\scripte_n$ by $\tilde{\scripte}_n$).
We also absorbed the error $\scripte_n$ coming from estimates of oscillatory integrals parallel to Lemma~\ref{aux1}
into $\tilde{\scripte}_n.$ 
In particular, 
\begin{equation}\label{norm3}
\|R(x,k_1^{(n)})\|_{L^2[x_{\pi-\alpha_1},(f_1)^{-1}g]}^2 = \frac{1}{2} f_1^{-1} \left( 1+ 
O(\tilde{\scripte}_n, g \tilde{g}^{-1}) \right)
\end{equation}
and, given the upper bound on $x_{\pi-\alpha_1}$ of Lemma~\ref{el2}, 
\begin{equation}\label{dec3}
R^2(f_1^{-1}g,k_1^{(n)}) \leq  e^{-g/2}(1+C(\tilde{\scripte}_n+ \tilde{g}^{-1} g)).
\end{equation}

4. Finally, we consider the solution on $[f_1^{-1}g, \infty).$ Here, 
\[ V_1(x) = -{\rm min}\left( \delta k_1,
\frac{g}{2x} \right) \frac{k_1^{(n-1)}}{4} (\sin 2\theta_1 + \sin 2\theta_2). \]
Equations for $\delta \theta_1$ and $R(x,k_1^{(n)})$ can be written in the following form:
\begin{eqnarray}\nonumber
(\log R^2(x, k^{(n)}_1 )'= -\frac{k_1^{(n-1)}}{8k_1^{(n)}} {\rm min}\left(\delta k_1, \frac{g}{2x}\right) 
\left( 1+ \cos 2\delta \theta_1 -\cos 4\theta_1 - \cos 2(\theta_1 +\theta_2) \right) + \\
\label{R2} +\frac{W(x)}{k_1^{(n)}} \sin 2 \theta_1  \\
\nonumber
(\delta \theta_1)'= \delta k_1 (1+O(\delta k_1)) + 
\frac{k_1^{(n-1)}}{8k_1^{(n)}} {\rm min}\left(\delta k_1, \frac{g}{2x}\right) \sin 2\delta \theta_1
(1-\cos 2(\theta_1+\theta_2))+ \\
\label{an2} + \frac{W(x)}{2k_1^{(n)}}  (\cos 2\theta_2 - \cos 2\theta_1). 
\end{eqnarray}
First of all, any contribution to the solution $R$ of the equation \eqref{R2} coming from the last term involving $W$ or two 
last summands in brackets in the first term is of the order $O(\scripte_n)$ by the same estimates as before
(in Lemma~\ref{aux1}). Therefore,
\begin{equation}\label{Rest2}
\frac{R(x,k_1^{(n)})^2}{R(f_1^{-1}g,k_1^{(n)})^2} = e^{- \frac18{\rm min}(\delta k_1, \frac{g}{2x}) 
\int_{f_1^{-1}g}^x (1+\cos 2\delta \theta_1)\,dy (1+O(\delta k_1))} (1+O(\scripte_n)).
\end{equation} 
Furthermore, same estimates show that for any $x>f_1^{-1}g,$ 
\begin{equation}\label{anest2}
\int\limits_{f_1^{-1}g}^x W(y) (\cos 2\theta_2 -\cos 2 \theta_1)\,dy = O(\scripte_n)
\end{equation}
as well. Consider a sequence of points $y_n,$ where $y_1 = f_1^{-1}g,$ and 
\[ y_n = {\rm min} \{ y: \delta \theta_1 (y) = \pi n + \delta \theta_1(y_1) \}. \]
According to \eqref{an2}, \eqref{anest2} and \eqref{genas},  we have 
$y_n - y_{n-1} \leq (2\pi +O(\scripte_n)) (\delta k_1)^{-1},$
and in each $(y_{n-1},y_n)$ there is an interval $I_n$ of length at least $(\pi/4 + O(\scripte_n))(\delta k_1)^{-1}$
where $\cos 2\delta \theta_1 >0.$ Let us assume for simplicity that the error terms in the above estimates do 
not consume more than half of the principal terms. 
Then from \eqref{R2} we see that on $(y_n, y_{n+1})$
\begin{equation}\label{decR1}
 R(y, k_1^{(n)})^2 \leq CR(y_1, k_1^{(n)})^2 e^{- \sum_{j=1}^{n-1} b_j}, 
\end{equation}
where 
\begin{equation}\label{bn}
 b_n = \frac{1}{16}\int\limits_{I_n} {\rm min}(\delta k_1, \frac{g}{2x}) \,dx \geq C
{\rm min} \left( 1, \frac{g}{g \tilde{g}^{-1}+ 3\pi n} \right).
\end{equation}
 Taking into account \eqref{dec3} and the estimates on $b_n,$ $y_n$ we have 
\begin{equation}\label{Rnbo11}
 \|R(y,k_1^{(n)})\|^2_{L^2[f_1^{-1}g, \infty)} \leq Ce^{-g/2} \sum\limits_{n=0}^\infty
e^{- \sum_{j=1}^n b_j } (\delta k_1)^{-1}. 
\end{equation}
Recall that $\delta k_1 = f_1^{-1} e^{g^{3/4}}.$ Provided that $g$ is large enough so that 
the sum in \eqref{Rnbo11} is finite, we can write
\begin{equation}\label{norm4}
\|R(y,k_1^{(n)})\|^2_{L^2[f_1^{-1}g, \infty)} = f_1^{-1} O(e^{-g/6}).
\end{equation}

Combining the estimates 
\eqref{xpi2norm}, \eqref{xpi2al}, \eqref{norm3} and \eqref{norm4}, we get the Splitting Lemma.
\end{proof}


\section{A Continuity Lemma} \label{continuity}

In this section, we prove an auxiliary technical lemma that we are going to need in the construction. 
The estimates here are fairly straightforward, if technical. Moreover, we choose simplicity over sophistication 
  and prove a version of the result which is sufficient for the application we have in mind.
With more extensive technical effort, the estimates in this section can be significantly improved.

Assume that the potential $V^{(n)}(x)$ is given on an interval $(x_n, \infty)$ by 
$V^{(n)}(x) 
= \sum_{j=1}^{2^{n-1}} V_j^{(n)}(x),$ where similarly to \eqref{sppo} we have 
\begin{equation}\label{sppo1}
V^{(n)}_j(x) = \left\{ \begin{array}{ll} -2 f_j^{(n)} k_j^{(n-1)} \sin(\theta^{(n)}_{2j-1}+\theta_{2j}^{(n)}),
 & x_n \leq x \leq x^{(n)}_{j,\pi-\alpha^{(n)}_j} \\
2 f_j^{(n)} k_j^{(n-1)} \sin(\theta^{(n)}_{2j-1}+\theta_{2j}^{(n)}), 
& x^{(n)}_{j,\pi-\alpha^{(n)}_j} \leq x \leq x_n+(f^{(n)}_j)^{-1}g_n \\
-{\rm min} \left( \delta \knj, \frac{g_n}{2x} \right) \frac{k_j^{(n-1)}}{4} (\sin 2\theta_{2j-1}^{(n)}
+\sin 2\theta_{2j}^{(n)}), &
x > x_n+(f_j^{(n)})^{-1}g_n. \end{array} \right.
\end{equation}
It is assumed that $\delta \theta^{(n)}_j(x_n) <\pi/2;$  $x^{(n)}_{j,\pi-\alpha^{(n)}_j}$ is then defined as
a minimal value of $x>x_n$ where $\delta \theta^{(n)}_j(x)$ reaches $\pi - \alpha^{(n)}_j.$
We adopt this definition in order to avoid making too many assumptions on what  happens for $x \leq x_n,$ 
although in the construction process that appears in the following sections it is straightforward to check that 
 $\delta \theta_j^{(n)}(x)$ does not reach $\pi -\alpha^{(n)}_j$ for $x \leq x_n.$ 
Similarly, in the actual construction we will always have $x_n << (f_j^{(n)})^{-1}g_n.$
As usual, $\delta k_j^{(n)} = \tilde{g}_{n}^{-1} f_j^{(n)}.$ For simplicity, we are also going to assume 
that 
\begin{equation}\label{stbop}
\|V^{(n)}\|_{L^\infty(\reals^+)} \leq \frac12 {\rm min}_j \{ (k_j^{(n)})^2 \}.
\end{equation} 
We also suppose that the condition \eqref{assplit} of the Splitting Lemma holds:
\begin{equation}\label{assplit1}
4\sum_j f^{(n)}_j k^{(n-1)}_j < \frac13 a_n,\,\,\,\,\, g_n >>1, \,\,\,2^n g_n a_n^{-2} 
\sum\limits_{j=1}^{2^{n-1}} f^{(n)}_j<<1,
\,\,\,\,\, \delta k^{(n)}_j < \frac{1}{12} a_n
\end{equation} 
for all $j.$ 
The potential $V^{(n)}(x)$ 
on $(0,x_n)$ comes from the previous steps, but for now 
we make no more specific assumptions about it. 
 Recall our notation 
$a_{n+1} = {\rm min}_{j,j'} |k_j^{(n)}-k_{j'}^{(n)}|.$
Let $\tilde{x}_n$ be such that for any $x \geq \tilde{x}_n,$ 
\begin{equation}\label{Vcon}
|V^{(n)}(x)| \leq \frac18 a_{n+1}{\rm min}_j \{ (\knj) \}.
\end{equation}
Also, take $\tilde{x}_n > g_n {\rm max}_j \{ (f_j^{(n)})^{-1}\}$ for convenience. 
\begin{lemma}\label{auxc1}
Assume that $V^{(n)}(x)$ is given by \eqref{sppo1} for $x>x_n.$
Then 
\begin{equation}\label{anbound2}
 \dthk (x,k^{(n)}_j) \geq P_n 
\left( \frac{x}{\tilde{x}_n} \right)^{\beta(\knj)} 
\end{equation}
for any $x \geq 2\tilde{x}_n.$ Here $\beta(k) = \frac{k_j^{(n-1)} g_n}{16k}.$
Moreover, if $k$ is such that 
\begin{equation}\label{dknotbig}
|k-k_j^{(n)}| \leq \frac{1}{4} a_{n+1},
\end{equation}
 then 
\begin{equation}\label{anbound1}
\dthk (x,k) \leq D_n \left( 1 + 
\left( \frac{x}{\tilde{x}_n} \right)^{\beta(k)} \right) 
\end{equation}
for any $x.$ 
Here $D_n,$ $P_n$ are positive constants 
defined only in terms of  $V^{(n)}(x)$ for $x \leq x_n,$ 
$x_n,$ $\{ f^{(n)}_j \},$ $a_{n+1}$ and $g_n$.
\end{lemma}
\begin{proof}
It is clear that the value of $\tilde{x}_n$ depends only on $x_n,$ $f^{(n)}_j,$ $a_{n+1}$ and $g_n,$
and can be easily estimated in terms of these quantities: 
\[ \tilde{x}_n \leq C{\rm max} \left(2^{n}g_n a^{-1}_{n+1}, x_n + g_n {\rm max}_j\{(f_j^{(n)})^{-1})\} \right). \] 
Notice that from \eqref{pruferan} it follows that 
\[ \left( \dthk(x,k) \right)' = 1 +\frac{1}{k^2} V(x) (\sin  \theta(x,k))^2 - \frac1k V(x) \sin 2\theta(x,k) \left( \dthk (x,k) \right). \]
Along with the boundary condition $\dthk(0,k)=0$ (which follows from the fact that we consider solutions 
satisfying fixed boundary condition for all $k$) this implies in our setting 
\begin{equation}\label{impan}
\dthk(x,k) = 
\int\limits_0^x  e^{-\frac1k \int_y^x V^{(n)}(s) \sin(2\theta(s,k))\,ds }
\left(1+ \frac{1}{k^2} V^{(n)}(y) (\sin \theta)^2 \right) \,dy.
\end{equation}
For any $k$ satisfying the assumption \eqref{dknotbig} of the lemma we have for $\tilde{x}_n \leq y \leq x$ 
by an estimate parallel to that in Lemma~\ref{aux1} and using \eqref{Vcon} and \eqref{assplit1}
\begin{equation}\label{keyint}
\int\limits_{y}^x V^{(n)}(y) \sin 2 \theta (y,k) \,dy = -\frac{k_j^{(n-1)} g_n}{16} \int\limits_y^x
\frac{ \cos(2\theta(s,k) - 2\tnj(s))}{s}\,ds +O(2^n a_{n+1}^{-2} g_n \sum_j f_j^{(n)}).
\end{equation} 
Equations \eqref{impan}, \eqref{keyint} together imply  a very rough estimate
\[ \dthk (x,k) \leq C \tilde{x}_n e^{\frac{1}{k} 
\int\limits_0^{\tilde{x}_n} |V^{(n)}(y)|\,dy} +
C e^{C2^n a^{-2}_{n+1} g_n \sum_j f_j^{(n)}} \int\limits_{\tilde{x}_n}^x e^{\beta(k) \log(x/y)}\,dy, \]  
where $C$ is a universal constant. Similarly, 
\[ \dthk (x, \knj) \geq C_1 e^{-C2^n a^{-2}_{n+1} g_n \sum_j f_j^{(n)}} \int\limits_{\tilde{x}_n}^x 
e^{\beta(\knj) \log(x/y)}\,dy. \]
Notice that 
\[ \int\limits_{\tilde{x}_n}^x 
e^{\beta(k) \log(x/y)}\,dy = \frac{\tilde{x}_n}{\beta(k)-1} \left( \left( \frac{x}{\tilde{x}_n} \right)^{\beta(k)}
 -\frac{x}{\tilde{x}_n} \right). \]
The parameter $\beta(k)$ is large because of $g_n,$ and we can assume for simplicity  $\beta(k)>2$ for all $k$ 
of interest (lying in a compact set where all $\knj$ lie). Then we get the result of the lemma. We can take, for example, 
\begin{eqnarray}\label{Dcon}
D_n = C \tilde{x}_n^2 e^{\frac{1}{k} \int\limits_0^{\tilde{x}_n} 
|V^{(n)}(y)|\,dy+ C2^n a^{-2}_{n+1} g_n \sum_j f_j^{(n)}}, \\ 
\label{Pcon}
P_n = C_1 \frac{\tilde{x}_n}{\beta(\knj)-1}  e^{-C2^n a^{-2}_{n+1} g_n \sum_j f_j^{(n)}}. 
\end{eqnarray}
\end{proof}

Now we prove the principal result of this section. 
\begin{lemma}[Continuity Lemma]\label{cont}
Choose $\tilde{g}_{n+1}$ so that $\tilde{g}_{n+1}^{1/4} \geq 2^n D_n P_n^{-1}.$ 
Then there exists $x_n'$ such that for any $x_{n+1} > x_n'$ the following holds. 
If  
\begin{equation}\label{kcloseco}
 |k- \knj| \leq 2 \tilde{g}_{n+1}^{-1/2} \left( \dthk (x_{n+1}, \knj) \right)^{-1}, 
\end{equation}
then 
\begin{equation}\label{Rclose}
R(x,k)^2 = R(x, \knj)^2 (1+O(\tilde{g}_{n+1}^{-1/4}))
\end{equation}
and 
\begin{equation}\label{tclose}
\dthk (x,k) = \dthk (x, \knj) (1+O(\tilde{g}_{n+1}^{-1/4}))
\end{equation}
for all $x \leq x_{n+1}.$ 
\end{lemma}
\begin{proof}
From \eqref{anbound2} it follows that $\dthk (x, k^{(n)}_j) \rightarrow \infty$ as $x \rightarrow \infty,$
and thus we can choose $x_n'$ so that for any $x_{n+1} \geq x_n'$ the values of $k$ satisfying \eqref{kcloseco}
also satisfy \eqref{dknotbig}. Then Lemma~\ref{auxc1} applies, and thus we have for any $x$ 
\[ \frac{\dthk (x,k) }{ \dthk (x_{n+1}, \knj)} \leq D_n P_n^{-1} \left( 1+ \left( \frac{x}{\tilde{x}_n} \right)^{\beta(k)}
\right) \left( \frac{x_{n+1}}{\tilde{x}_n}\right)^{-\beta (\knj)}. \]
Hence, by definition of $\tilde{g}_{n+1},$ we have   
\begin{equation}\label{andiff11}
 |\theta(x,k) - \theta (x, \knj)| \leq 2^{-n} \tilde{g}_{n+1}^{-1/4} \left( 1+ \left( \frac{x}{\tilde{x}_n} \right)^{\beta(k)}
\right) \left( \frac{x_{n+1}}{\tilde{x}_n}\right)^{-\beta (\knj)}. 
\end{equation}
Recall from \eqref{pruferam} that
\[ R(x,k)^2/R(x,\knj)^2 = e^{\frac{1}{k} \int\limits_0^x V^{(n)}(y) \sin 2\theta(y,k) \,dy - 
\frac{1}{\knj} \int\limits_0^x V^{(n)}(y) \sin 2\theta(y,\knj) \,dy}. \]
First of all, we can replace $1/k$ with $1/\knj$ in the first term in the exponent at the expense of creating an error factor 
of size at most 
\[ e^{C |k-\knj| \int\limits_0^{x_{n+1}} |V^{(n)}(y)|\,dy} \]
for $x \leq x_{n+1}.$ From \eqref{kcloseco} and a bound \eqref{anbound2} of Lemma~\ref{auxc1} it is clear that we can choose $x_n'$ so that 
this factor does not exceed $1+O(\tilde{g}_{n+1}^{-1/4}).$ 
Next, we need to bound
\begin{equation}\label{finerr3}
 \frac{1}{\knj} 
\int\limits_0^x V^{(n)}(y) |\theta (y,k) -\theta(y,\knj)| \,dy. 
\end{equation} 
Using \eqref{andiff11}, we can estimate the expression \eqref{finerr3} for $x \leq x_{n+1}$ by 
\begin{equation}\label{simp12}
C 2^{-n} \tilde{g}_{n+1}^{-1/4} x_{n+1}^{-\beta(\knj)} \tilde{x}_n^{\beta(\knj)-\beta(k)}
\int\limits_0^{x_{n+1}} y^{\beta(k)} |V^{(n)}(y)| \,dy. 
\end{equation}
The contribution to \eqref{simp12} 
coming from the integral over $[0,\tilde{x}_n]$ can be made $O(\tilde{g}_{n+1}^{-1/4})$ by the
choice of $x_n'.$
Using an inequality $|V^{(n)}(y)| \leq C 2^n \beta(k) y^{-1}$ for $y> \tilde{x}_n,$ we estimate the remaining part of \eqref{simp12}
by $C \tilde{g}_{n+1}^{-1/4} (x_{n+1}/\tilde{x}_n)^{\beta(k)-\beta(\knj)}.$   
Notice that from the definition of $\beta(k),$ we have 
\begin{equation}\label{betabo}
|\beta (k)- \beta (\knj)| \leq Cg_n  |k-\knj| \leq C g_n \tilde{g}_{n+1}^{-1/2} \left( \dthk (x_{n+1}, \knj)\right)^{-1}.
\end{equation}
By \eqref{anbound2} of Lemma~\ref{auxc1} and the definition of $\tilde{g}_{n+1},$ the right hand side in \eqref{betabo}
does not exceed $C (x_{n+1}/\tilde{x}_n)^{-\beta(\knj)}.$ Thus,
\begin{equation}\label{finest11}
\left(\frac{x_{n+1}}{\tilde{x}_n} \right)^{\beta(k)-\beta(\knj)} 
\leq e^{C (\log (x_{n+1}/\tilde{x}_n)) (x_{n+1}/\tilde{x}_n)^{-\beta(\knj)}},
\end{equation}
which is bounded by some fixed universal constant.
Combining the estimates, we see that if we choose $x_n'$ to satisfy all of the above requirements, then 
\eqref{Rclose} holds for any $x \leq x_{n+1},$ provided that $x_{n+1} \geq x_n'.$

Notice also that from \eqref{andiff11} and \eqref{finest11}, it follows that 
for $x \in (0,x_{n+1})$ we have $|\theta(x,k) - \theta(x, \knj)| \leq C \tilde{g}_{n+1}^{-1/4}.$   
Thus \eqref{tclose} follows from \eqref{Rclose} and \eqref{impan}, which can also be written as
\begin{equation}\label{animp1}
\dthk (x,k) = R(x,k)^{-2} \int\limits_0^x R(y,k)^2 \left( 1+ \frac{1}{k^2} 
V^{(n)}(y) (\sin \theta)^2 \right) \,dy.
\end{equation}
We used \eqref{stbop} in the last estimate.
\end{proof}

\section{A Brick of Construction} \label{brick}

Our goal is  to organize a Cantor-like process, splitting each eigenvalue we have on the $n$th step 
into  two on the $(n+1)$th step. To do this we need to see how norm splitting works in a more general setting than 
Splitting Lemma, which starts at the origin. The needed estimates are provided by the following lemma, which helps 
to connect different steps.

Assume that for some $j,$ we have for $k$ satisfying 
\[ |k- k_j^{(n)}| \leq 2\tilde{g}_{n+1}^{-1/2}   \left( \dthk (x_{n+1}, k_j^{(n)})\right)^{-1} \]
the bounds of Continuity Lemma:
\begin{equation}\label{cR}
 R(x,k)^2 = R(x,k_j^{(n)})^2 (1+ O(\tilde{g}^{-1/4}_{n+1})) 
\end{equation}
and 
\begin{equation}\label{can}
 \dthk (x,k) = \dthk (x, k_j^{(n)})(1+O(\tilde{g}_{n+1}^{-1/4})) 
\end{equation}
for all $x \leq x_{n+1}.$ 
Set 
\begin{eqnarray}
\label{fs}  f_j^{(n+1)} = \tilde{g}^{1/2}_{n+1} \left( \dthk (x_{n+1}, k_j^{(n)})\right)^{-1}  \\
\label{dks} \delta k_j^{(n+1)} = \tilde{g}^{-1/2}_{n+1} \left( \dthk (x_{n+1}, k_j^{(n)})\right)^{-1} ,
\end{eqnarray}
and let $k_{2j}^{(n+1)},$ $k_{2j-1}^{(n+1)}$ be the ends of the interval of size $\delta k_j^{(n+1)}$ 
centered at $k^{(n)}_j.$ 
For $x>x_{n+1},$ define potential $V^{(n+1)}(x)$ to be given by \eqref{sppo1} with $n+1$ replacing $n.$
We are going to assume also that the conditions \eqref{assplit1} of the Splitting Lemma are satisfied
(with $n$ replaced by  $n+1$ in our current situation). 
Define a parameter $A_j^{(n)}$ by 
\begin{equation}\label{Adef} 
A_j^{(n)}= \dthk (x_{n+1}, k_j^{(n)}) R(x_{n+1}, k_j^{(n)})^2.
\end{equation}
We have 
\begin{lemma}[Connection Lemma]\label{connect}
Under the above assumptions, $k_{2j-1}^{(n+1)}$ and $k_{2j}^{(n+1)}$ are eigenvalues of $H_{V^{(n+1)}},$ 
and the norms of the corresponding eigenfunctions on $(x_{n+1}, \infty)$ satisfy 
\begin{equation}\label{norcrit}
\|R(x,k^{(n+1)}_{2j,2j-1})\|^2_{L^2(x_{n+1}, \infty)} = A_j^{(n)}(1+ O(\tilde{\scripte}_{n+1}, \tilde{g}^{-1/4}_{n+1})).
\end{equation}
\end{lemma}
\begin{proof}
To simplify notation, let us assume without loss of generality that $j=1.$
According to our assumptions, $\delta \theta_1^{(n+1)} (x_{n+1}) = \tilde{g}_{n+1}^{-1/2}(1+O(\tilde{g}_{n+1}^{-1/4})).$ 
If $\delta \theta_1^{(n+1)} (x_{n+1})$ were zero, we would be exactly in a situation of the Splitting Lemma
construction, and would have two eigenfunctions of norm $\sim 2(f_j^{(n+1)})^{-1} R(x_{n+1}, k_j^{(n)})^2.$
We have to take into account the fact that on the interval where $\delta \theta_1^{(n+1)}(x)  \in [0,\delta \theta_1^{(n+1)} (x_{n+1})]$
 a significant decay takes place in the Splitting Lemma, and we lose this decay in our situation. 

Consider an auxiliary problem for $\overline{R},$ $\delta \overline{\theta}$ satisfying same equations \eqref{keyR}, 
\eqref{keyan} as 
$R(x,k_1^{(n+1)}),$ $\delta \theta^{(n+1)}_1$ but with 
$\overline{V}^{(n+1)}$ given by \eqref{sppo} with $n+1$ instead of $n.$ In other words, consider the case of $(n+1)$st step potential 
starting right away at the origin.
The Small Angle and Splitting Lemmas imply that 
\begin{equation}\label{auxss}
\|\overline{R}(x,k_1^{(n+1)})\|^2_{L^2(x_{\delta \theta_1(x_{n+1})},  \infty)} = (f_1^{(n+1)})^{-1}(1+ O(\tilde{\scripte}_{n+1},
\tilde{g}_{n+1}^{-1/2} g_{n+1} ) ). 
\end{equation}
But the same estimates as in these Lemmas apply to $R(x,k^{(n+1)}_1),$ $\delta \theta_1$
 on $(x_{n+1}, \infty),$ by identical arguments. 
The only  adjustment we need to make is by a factor 
\[ R(x_{n+1}, k_1^{(n+1)})^2/\overline{R}(x_{\delta \theta_1 (x_{n+1})},k_1^{(n+1)})^2. \]
By Lemma~\ref{sman} (Small Angle Lemma),
\begin{equation}\label{dloss3}
 \overline{R}(x_{\delta \theta_1 (x_{n+1})},k_1^{(n+1)})^2=
 \left( \frac{\delta k^{(n+1)}_1 }{\delta \theta_1^{(n+1)}(x_{n+1}) f_1^{(n+1)}+\delta k_1^{(n+1)}} \right) 
\left( 1+ O(\scripte_{n+1}, 
\tilde{g}_{n+1}^{-1}g_{n+1})\right). 
\end{equation}
According to \eqref{auxss}, \eqref{dloss3} we have 
\begin{equation}\label{normsp1}
\|R(x,k_{1,2}^{(n+1)})\|^2_{L^2(x_{n+1},\infty)} =\frac{ R(x_{n+1}, k_1^{(n)})^2}{f_1^{(n+1)}} 
\left( \frac{\delta \theta_1^{(n+1)}(x_{n+1}) f_1^{(n+1)}}{\delta k^{(n+1)}_1 }+1 \right)
\left( 1+O(\tilde{\scripte}_{n+1}, \tilde{g}_{n+1}^{-1/4}) \right). 
\end{equation}
We utilized assumption \eqref{cR} in the error estimate, and also used $\tilde{g}_{n+1}^{-1}g_{n+1} <
C\tilde{g}_{n+1}^{-1/4}.$ 
Because of the assumption \eqref{can}, we have 
\begin{equation}\label{add33}
  \frac{\delta \theta_1^{(n+1)}(x_{n+1})}{\delta k_1^{(n+1)}} = \dthk (x_{n+1}, k_1^{(n)})(1+O(\tilde{g}_{n+1}^{-1/4})). 
\end{equation}
The term $(f_1^{(n+1)})^{-1} R(x_{n+1}, k_1^{(n)})^2 = A_j^{(n)} \tilde{g}^{-1/2}_{n+1},$ 
and hence contributes a correction of higher order.
Thus by definition of $A_j^{(n)}$ and \eqref{add33}, \eqref{normsp1} leads exactly to \eqref{norcrit}.
\end{proof}

The usefulness of Lemma~\ref{connect} becomes clear if one looks at \eqref{animp1}, which in particular implies
\[ \dthk (x,k) = R(x,k)^{-2} \|R(y,k)\|^2_{L^2(0,x)}(1+O(\|V\|_{L^\infty})). \]
Thus in the asymptotic regime where $V$ is small, one can expect, informally, 
\[ \dthk (x,k) \sim R(x,k)^{-2} \|R(y,k)\|^2_{L^2(0,x)},\]
and the role of $A$ is played by $\|R(x,k)\|^2_{L^2(0,x)}.$ Lemma~\ref{connect} allows to get two eigenvalues with the 
norm on $(x,\infty)$ same as $A=\|R(x,k)\|^2_{L^2(0,x)}.$ This will provide the mechanism for norm doubling in the 
asymptotic regime which is crucial for the construction.  
 
\section{The Construction} \label{construction}

Now we are ready to provide a precise description of the construction of our potential $V(x).$
In the beginning, we take a potential given by \eqref{sppo},
like in Splitting Lemma with $n=1$.   
Start by choosing $\delta k_1^{(1)},$ $k^{(0)}_1,$ $g_1$ and $f^{(1)}_1.$
The precise choice of the parameters
on the first step of construction is not so important, but we are going to keep  it similar to the 
subsequent steps.  
 To satisfy \eqref{genas}
we choose $k^{(0)}_1$ sufficiently large and $\delta k^{(1)}_1$ small compared to $k^{(0)}_1.$ 
We take $g_1$ large,
impose a relationship $\delta k_1^{(1)} = \tilde{g}^{-1}_1 f_1^{(1)}$ and 
take $e^{g_1}f_1^{(1)}$ small. The smallness of the latter number and $\delta k_1^{(1)}$ can be achieved
by taking $f_1^{(1)}$ sufficiently small while keeping $g_1$ fixed.
In particular, we take $g_1$ large enough and $f_1^{(1)}$ small enough for the conditions \eqref{assplit1}
of the Splitting to be satisfied.  
The only additional restriction we need to impose on the 
parameters comes from the condition $|V(x)| \leq h(x)/(1+x).$ 
By the monotonicity of $h$ and the definition of potential, the inequality 
\begin{equation}\label{ubp33}
 2f_1^{(1)}k^{(0)}_1 \leq \frac{h((f_1^{(1)})^{-1} g_1)}{1+(f_1^{(1)})^{-1} g_1}, 
\end{equation}
implies the needed bound for $x \geq (f_1^{(1)})^{-1} g_1.$ It can be rewritten as
 $h((f_1^{(1)})^{-1} g_1) \geq 2 (g_1+f_1^{(1)}) k^{(0)}_1$ which holds for sufficiently small $f_1^{(1)}.$ 
To simplify consideration for smaller $x,$ we are going to assume without loss of generality that 
$h(x)/(1+x)$ is decreasing. Then \eqref{ubp33} also implies the bound for $x< (f_1^{(1)})^{-1} g_1.$ 
The Splitting Lemma then gives us two eigenvalues at $k^{(1)}_{1,2}$ with corresponding eigenfunctions 
having norms $(f_1^{(1)})^{-1} (1+O(e^{g_1}f_1^{(1)}, \tilde{g}_1^{-1/2}g_1)).$

Let us now describe the induction process. We are going to assume that on the $n$th step we have the potential
$V^{(n)}(x),$ in particular given by \eqref{sppo1} on $(x_n, \infty).$ We assume that this potential leads to eigenvalues
$k^{(n)}_j,$ $j =1,\dots 2^n$ and the following holds: \\
1. The norms of the corresponding eigenfunctions satisfy for every $j$
\begin{equation}\label{keyef}
\| R(x,k_j^{(n)}) \|^2_{L^2(0,\infty)} \leq C_1 2^n\prod\limits_{l=1}^n (1+ Er_l),
\end{equation} 
where $|Er_l| \leq C_2 2^{-l},$ and $C_{1,2}$ are fixed constants which do not depend on $n,l.$  \\
2. The angle derivative in energy satisfies for every $j$ 
\begin{equation}\label{dankey}
\dthk (x, \knj) \leq C_1 2^n R(x, \knj)^{-2} \prod\limits_{l=1}^n (1+ Er_l)
\end{equation}
for any $x.$ \\
3. The conditions \eqref{assplit1} of the Splitting Lemma are satisfied, and $|V^{(n)}(x)| \leq h(x)/(1+x).$ \\
4. We assume that \eqref{genas} holds; this will be nearly automatic because of the rapid 
decrease of $\delta k_j^{(n)}.$  \\

\begin{figure}
\begin{center}
\scalebox{0.75}{\includegraphics{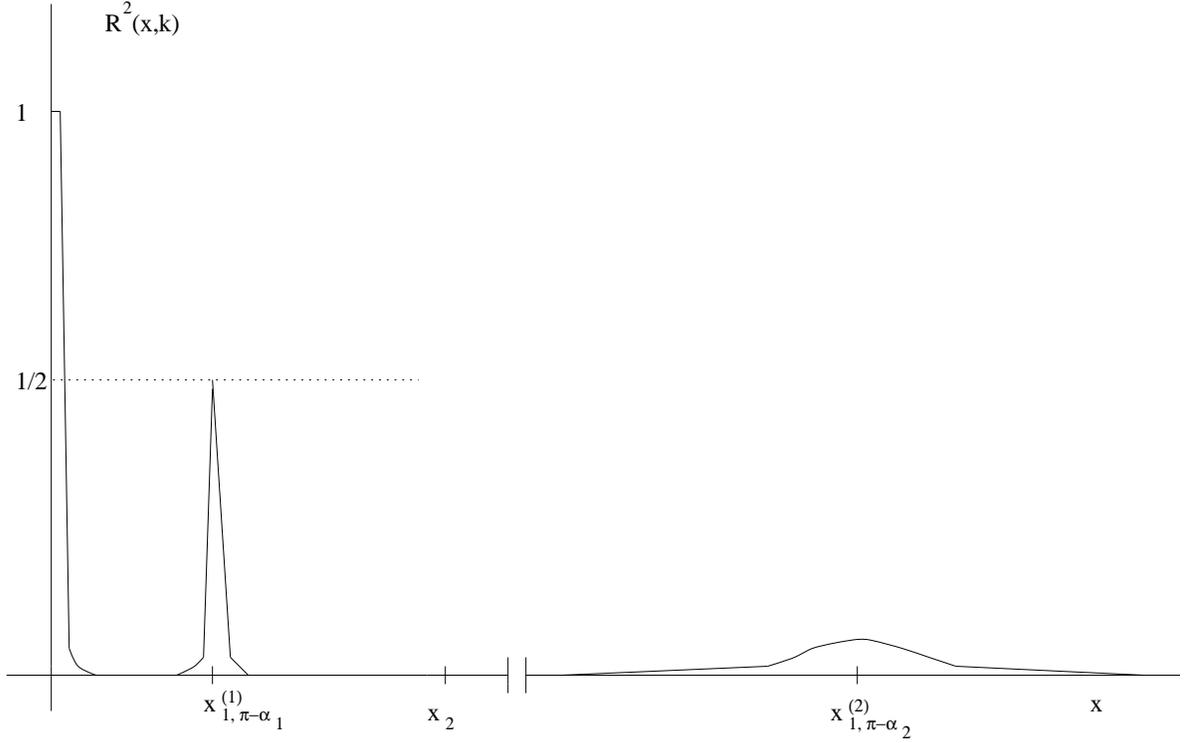}}
\caption{The structure of the eigenfunction at $k=k^{(n)}_1$. 
Only the first two stages are shown. The factor by which the eigenfunction 
grows from the splitting point $x_n$ to its top value at $x^{(n)}_{1,\pi-
\alpha_1^{(n)}}$ is 
$\sim (1/2) \tilde{g}_{n}^{1/2}.$ This can be seen from the argument in
Connection Lemma. Each bump contains roughly as much norm as everything 
before it.}
\label{struct1}
\end{center}
\end{figure}

Clearly for $n=1$ these assumptions are satisfied with $x_n=0$ and sufficiently large $C_{1,2}.$ $C_1,$ 
in particular, should be taken $\sim (f_1^{(1)})^{-1}$ or larger. The fact that \eqref{dankey} is satisfied for
$n=1$ follows from \eqref{keyef} and \eqref{animp1}. 
The structure of the eigenfunctions that we are constructing 
is illustrated by Figure~\ref{struct1}.

Now let us make an induction step. 
Choose $\tilde{g}_{n+1}$ so 
that $\tilde{g}_{n+1}^{1/4} > 2^n  D_n P_n^{-1},$ where 
$D_n,P_n$ are given by  \eqref{Dcon}, \eqref{Pcon}. We should also make sure that $g_{n+1}$ is large
enough so that the second condition of \eqref{assplit1} 
needed for Splitting and Connection Lemmas holds; but this is basically automatic after the first step, 
since  $g_n$ is virtually forced by construction to be an increasing sequence. 
We require in addition that $\tilde{g}_{n+1}^{-1/4} \leq C2^{-n-1}$ (which is again basically automatic).   
For each $j=1,\dots, 2^n,$ choose 
\begin{equation}\label{parch}
f^{(n+1)}_j = \tilde{g}_{n+1}^{1/2} \left( \dthk (x_{n+1}, \knj ) \right)^{-1}, \,\,\,\,\,
\delta k^{(n+1)}_j = \tilde{g}_{n+1}^{-1/2} \left( \dthk (x_{n+1}, \knj ) \right)^{-1}, 
\end{equation}
where the splitting point $x_{n+1}$ remains the only parameter of the $(n+1)$st step to be chosen. 
Define 
\begin{equation}\label{tf33}
 \tilde{f}_{n+1} = Ca_{n+1}^2 2^{-3(n+1)} e^{-g_{n+1}}. 
\end{equation}
One condition that we impose on $x_{n+1}$ is that $f^{(n+1)}_j$ defined by \eqref{parch}
are smaller than $\tilde{f}_{n+1}$ for every $j.$ This can always be achieved by taking $x_{n+1}$ large
enough because of \eqref{anbound2}. Notice that such choice ensures that $\tilde{\scripte}_{n+1}
\leq C 2^{-n-1}.$ Also, the constant $C$ is \eqref{tf33} is chosen so that the first, third and fourth
conditions of \eqref{assplit1} needed for Splitting and Connection Lemmas on the $(n+1)$st step hold.

The second condition imposed on $x_{n+1}$ is that Lemma~\ref{cont} (Continuity Lemma) applies,
that is, $x_{n+1} >x_n',$ $x_n'$ is defined in Lemma~\ref{cont}, and $\delta k_j^{(n)} \leq \frac{1}{12} a_{n+1}$ 
for all $j.$  


The third condition comes from the constraint $|V^{(n+1)}(x)| \leq h(x)/(1+x).$  
 The potential $V^{(n+1)}$ on $[x_{n+1},\infty)$ is going to be defined by \eqref{sppo1}
with $n+1$ instead of $n,$ so the latter requirement is going to be satisfied if for every $j,$
\begin{equation}\label{grreq}
f_j^{(n+1)} \leq 2^{-n-1} \frac{h(x_{n+1} + (f_j^{(n+1)})^{-1} g_{n+1})}
{x_{n+1} + (f_j^{(n+1)})^{-1} g_{n+1}}.
\end{equation}
By \eqref{anbound2} we have 
\[ \left( \dthk (x, \knj) /x \right) \stackrel{x \rightarrow \infty}{\longrightarrow} \infty \]
(since we assume that $g_n$ is large enough). Thus, for large enough $x_{n+1}$ 
 we have $x_{n+1} \leq (f_j^{(n+1)})^{-1} g_{n+1}$ (recall that $f_j^{(n+1)}$ is defined by 
\eqref{parch}). Then the condition \eqref{grreq} reduces to 
\[ 2^{n+2} g_{n+1} \leq  h(x_{n+1} + (f_j^{(n+1)})^{-1} g_{n+1}). \]
This can always be achieved by taking $x_{n+1}$ sufficiently large while $g_{n+1}$ is fixed.

By the choice of $\delta k_j^{(n+1)},$ $x_{n+1}$ and $\tilde{g}_{n+1},$ the result of  
Lemma~\ref{cont} (Continuity Lemma) applies to $k_j^{(n+1)},$ 
and  therefore we can use Lemma~\ref{connect} 
(Connection Lemma). 
Together, Continuity and Connection Lemmas show that $k_j^{(n+1)}$ are eigenvalues with 
eigenfunctions satisfying 
\[ \|R(x,k_j^{(n+1)})\|^2_{L^2} \leq C_1 2^n \prod\limits_{l=1}^{n} (1+Er_l) \left( (1+ C \tilde{g}_{n+1}^{-1/4})
+ (1+ C (\tilde{g}_{n+1}^{-1/4}+ \tilde{\scripte}_{n+1})) \right). \]
Here the first summand in the brackets comes from $(0, x_{n+1})$ and was estimated by the induction assumption 
\eqref{keyef} and Continuity Lemma, while the second summand comes from $(x_{n+1}, \infty)$ and was 
estimated using \eqref{dankey}, Continuity and Connection Lemmas. 
Altogether, we obtain 
\begin{equation}\label{finnorm3}
 \|R(x,k_j^{(n+1)})\|^2_{L^2} \leq C_1 2^{n+1} \prod\limits_{l=1}^{n+1} (1+Er_l),
\end{equation}
where $Er_{n+1} \leq C_2 2^{-n-1}$ and thus satisfies the required bound. 
Next,
\[ \dthk (x, k_j^{(n+1)}) = \dthk (x_{n+1}, k_j^{(n+1)}) \frac{R(x_{n+1},k_j^{(n+1)})^2}{R(x,k_j^{(n+1)})^2} +
 \int\limits_{x_{n+1}}^x \frac{R(y, k_j^{(n+1)})^2}{R(x,k_j^{(n+1)})^2} (1+O(|V(y)|))\,dy. \]
By induction assumption and Continuity Lemma, the first summand does not exceed
\[ R(x,k_j^{(n+1)})^{-2} C_1 2^n \prod\limits_{l=1}^n (1+ Er_l) (1+C \tilde{g}_{n+1}^{-1/4}). \]
The second summand, by Connection Lemma, does not exceed
\[ R(x,k_j^{(n+1)})^{-2} C_1 2^n \prod\limits_{l=1}^n (1+ Er_l) (1+C (\tilde{\scripte}_{n+1}+\tilde{g}_{n+1}^{-1/4}))
(1+C2^n \tilde{f}_{n+1}).  \]
Therefore, taking into account the above estimates on $\tilde{\scripte}_{n+1},$  $\tilde{g}_{n+1}^{-1/4}$ and the choice 
of $\tilde{f}_{n+1},$ we get that 
\[ \dthk (x, k_j^{(n+1)})= R(x,k_j^{(n+1)})^{-2} C_1 2^{n+1} \prod\limits_{l=1}^{n+1} (1+ Er_l), \]
with $E_{n+1} \leq C_2 2^{-n-1}.$ The third component of the induction step was already handled by the choice of $x_{n+1}.$ 
The fourth component is easily satisfied, since by inspection of the argument 
it is clear that $\delta k_j^{(n)}$ decay rapidly in $n.$ 
In fact we can assume freely that $g_n >2^n$ (there is no restriction on how large
$g_n$ we can choose on each step, only a bound from below). Since $f_j^{(n)} \leq C 2^{-3n} e^{-g_n}$ for any $j$ and $\delta k_j^{(n)}$  
is even smaller, we can clearly arrange for \eqref{genas} to be valid. 
This concludes the induction step. 

Notice that by construction, ${\rm max}_j \delta k^{(n+1)}_j \leq \tilde{f}_{n+1} \tilde{g}_{n+1}^{-1/2},$ and by definition 
 of $a_{n+1},$ $\tilde{f}_{n+1} \leq 2^{-3n} e^{-g_n} ({\rm min}_j  \delta k_j^{(n)})^2.$ It is easy to see from  these relations that 
for all sufficiently large $n,$ 
\begin{equation}\label{inter}
\sum\limits_{m=n+1}^\infty {\rm max}_j (\delta k_j^{(m)}) \leq \frac13 {\rm min}_j (\delta k_j^{(n)}). 
\end{equation}
Denote $I_j^{(n)}$ the intervals with centers at $\knj$ of size $\frac13 \delta \knj.$ Then by \eqref{inter}, for sufficiently large $n,$ 
$2I^{(n)}_j$ are disjoint for different $j$ (where by $2I^{(n)}_j$ we mean the intervals  centered at $\knj$ twice the size of $I^{(n)}_j$). 
At the same time, all eigenvalues $k_l^{(m)}$ that are generated from $\knj$ at a later steps and their intervals 
$I_l^{(m)}$ lie in $I^{(n)}_j$
(these are the eigenvalues $k^{(m)}_l,$ $l = (j-1)2^{m-n}+1,\dots j 2^{m-n}$ for any $m>n$).   
For the size of $I_j^{(n)}$ we have an estimate $|I_j^{(n)}| \leq Ce^{-2^n}$ provided that we choose $g_n >2^n$ (the estimate is even stronger,
but this will suffice). 

Let us summarize our findings.
\begin{theorem}\label{evthm}
Applying the construction described above we obtain a sequence of potentials $V^{(n)}(x)$ with the following properties. 
Each $V^{(n)}(x)$ satisfies $|V^{(n)}(x)| \leq h(x)/(1+x).$ If $m>n,$ then $V^{(m)}(x) = V^{(n)}(x)$ for $x \leq x_n,$ where 
$\{x_n\}$ is a strictly  increasing sequence  tending to infinity as $n$ grows. For each $n,$ the operator $H_{V^{(n)}}$ 
has $2^n$ eigenvalues $\knj,$ $j =1,\dots 2^n$ such that  for any $j$ 
\begin{equation}\label{cruc}
\int\limits_0^\infty R(x, k_j^{(n)})^2 \,dx \leq B^{-1} 2^{n},
\end{equation}
where $B>0$ is a universal constant independent of $n.$ 
Moreover, for each $\knj$ we can define an interval $I_j^{(n)}$ centered at $\knj$ such that $|I_j^{(n)}| \leq C e^{-2^n}$ 
and for any $m>n,$ all $k^{(m)}_l$ with $l  = (j-1)2^{m-n}+1,\dots j 2^{m-n}$ belong to $I_j^{(n)}.$  
At the same time, intervals $2I^{(n)}_j,$ $j=1, \dots, 2^n$ for fixed $n$ are disjoint from each other.
\end{theorem}
\begin{proof}
We have already shown every aspect of the theorem except we have not directly stated \eqref{cruc}. However this estimate
follows directly from the first induction assumption, \eqref{keyef}, since $Er_l \leq C_2 2^{-l}.$  
\end{proof}

\section{Singular Continuous Spectrum} \label{sc}

In this section we prove Theorem~\ref{main}. We are going to define potential $V(x)$ by a condition $V(x)=V^{(n)}(x)$ for $x \leq x_n.$ 
By Theorem~\ref{evthm}, this defines unambiguously $V(x)$ on all positive semi-axis. We prove the following result, which implies 
Theorem~\ref{main}.

\begin{theorem}\label{scsp}  
The singular continuous spectrum of the operator $H_V$ is nonempty. 
\end{theorem}
\begin{proof}
Let us denote the spectral measure of $H_V$ by $\mu,$ and spectral measures of $H_{V^{(n)}}$ by $\mu_n.$ 
Notice that since $V=V^{(n)}$ for $x <x_n$ and $x_n \rightarrow \infty,$ the operator $H_{V^{(n)}}$ converges
to $H_V$ in a strong resolvent sense. The strong resolvent convergence implies weak convergence of $\mu_n$ to $\mu$ 
(see, e.g. \cite{RS1}, Theorem VII.20). Consider any interval $I_j^{(n)}$ with $n$ sufficiently large. Then by Theorem~\ref{evthm} 
and Lemma~\ref{el}, we have 
$\mu_m (I_j^{(n)}) \geq B 2^{-n}$ for each $m>n.$ This implies by weak convergence that $\mu(2I_j^{(n)}) \geq B 2^{-n}$ 
for every $n$ and $j =1,\dots, 2^n.$  Assume that $\mu$ does not have singular continuous component, so that 
$\mu =  \mu^{{\rm ac}}  + \mu^{{\rm pp}}.$ Take a small number $\epsilon < B/10,$ then there exists a $\delta_\epsilon$ such that 
if $|S| < \delta_\epsilon,$  $\mu^{{\rm ac}}(S) < \epsilon$ (we denote here by $|S|$ the Lebesgue measure of the set  $S$). 
Choose $n_1$ so that for any $n>n_1,$ 
\[ \sum\limits_{j=1}^{2^n} |2I_j^{(n)}| \leq C2^{n+1} e^{-2^n} < \delta_\epsilon. \]
Then 
\begin{equation}\label{ac}
 \mu^{{\rm ac}}(\cup_{j=1} (2I_j^{(n)})) \leq \epsilon < B/10
\end{equation}
for any $n>n_1.$ Let $\mu^{{\rm pp}} =  \sum_{l=1}^\infty  \rho_l \delta(x-x_l).$  
Pick $N_\epsilon$ so that $\sum_{l=N_\epsilon}^\infty \rho_l < \epsilon.$ Choose $n_2$ so that $2^{n_2}>N_\epsilon.$ 
Then for any $n > n_2,$ there are at least $2^n -2^{n_2}$ intervals $2I^{(n)}_j$ which do not contain any of the $x_l$ 
with $l \leq N_\epsilon.$ Denote the set of such intervals $Q_n.$ Then we have 
\begin{equation}\label{pp}
\mu^{{\rm pp}}(\cup_{j \in Q_n} (2I_j^{(n)})) \leq \epsilon < B/10
\end{equation}
for any $n>n_2.$ 
At the same time, by Theorem~\ref{evthm}, we have
\[ \mu(\cup_{j \in Q_n} (2I_j^{(n)})) \geq B/2 \]
if $n>n_2.$ 
For $n \geq {\rm max}(n_1,n_2),$ this gives us a contradiction with \eqref{ac}, \eqref{pp} and our assumption on the absence 
of singular continuous component. 
\end{proof}

Theorem~\ref{main} follows from Theorem~\ref{scsp} immediately since from Theorem~\ref{evthm} we know that 
$|V(x)| \leq h(x)/(1+x).$ 

\section{Asymptotic Incompleteness of Wave Operators} \label{iwo}

In \cite{CK3}, it was proved that the M\"oller wave operators exist if $V \in L^p,$ $p<2,$ and in addition $\lim_{x \rightarrow \infty}
\int_0^x V(s)\,ds$ exists. Therefore, to prove Theorem~\ref{asin}, we only need to show that the potential $V(x),$
constructed in the previous sections, can be chosen to be conditionally integrable. 

\begin{proof}[Proof of Theorem~\ref{asin}]
Consider one of the components of $V^{(n)}(x),$ say $V_j^{(n)}(x),$ on $(x_n,\infty).$ 
Simple integration by parts similar to Lemma~\ref{aux1} shows that for any $x,$
\[ \left|\, \int\limits_{x_n}^x V_j^{(n)}(y) \,dy  \right| \leq C\left(f_j^{(n)} + \int\limits_{x_n}^\infty \tilde{f}_j^{(n)}(y) |V^{(n)}(y)|\,dy,\right) \]
where similarly to what we had before in \eqref{auxf}, $\tilde{f}_j^{(n)}$ is given by 
 \[ \tilde{f}_j^{(n)} (y) = \left\{ \begin{array}{ll} -f^{(n)}_j, & x_n \leq x \leq x^{(n)}_{j,\pi-\alpha^{(n)}_j} \\
f_j^{(n)}, & x^{(n)}_{j,\pi-\alpha_j^{(n)}} \leq x \leq x_n+(f_j^{(n)})^{-1}g_n \\
-{\rm min}(\delta k^{(n)}_j, \frac{g_n}{2x} ) \frac{k_j^{(n-1)}}{4}, & x > x_n+(f_j^{(n)})^{-1}g_n. 
\end{array}
\right. \]
Since we also have $|V^{(n)}(y)| \leq C \sum_j \tilde{f}_j^{(n)}(y),$ a simple computation gives that 
\[ \left| \,\int\limits_{x_n}^x V^{(n)}(y) \,dy \right| \leq C 2^n g_n \sum_j f_j^{(n)} \]
for any $x>x_n.$  
According to our choice of parameters in the construction, in particular $\tilde{f}_n,$ the bound on the right hand side does
not exceed $C 2^{-n}.$ Summing up the contributions from $(x_n,x_{n+1})$ for different $n,$ we get the result.
\end{proof}

\section{Discussion} \label{disc}

In this paper, we did not try to study fine properties of the imbedded singular continuous measure, such as Hausdorff dimension. 
The reasonable conjecture is that for $|V(x)| \leq C(1+x)^{-\alpha},$ there could be imbedded singular continuous  measures of dimension 
up  to $1-2\alpha$ (see \cite{Rem1,Re2}). One might try to use ideas of this paper to construct such measures. 
This, however, would have to lead to a different plan of attack and more technical effort. In particular, one should use 
approximations of the limiting measure by absolutely continuous, rather than by pure point, measures, cutting $V^{(n)}(x)$ off at an
appropriate point. In addition, one would have to be more careful with the separation of scales on different steps of construction, 
which is extensively used in our argument to control the errors. 
Such an approach would require an 
improvement in control of many estimates, in particular in continuity estimates of Section~\ref{continuity}. 
To keep the size of the paper in check, we did not try to pursue this direction here.

Another interesting open question is what happens sharply at the Coulomb rate of decay $1/x.$ In \cite{Rem1}, it is shown that imbedded singular
spectrum cannot occur if $V(x) = o(1/x)$ or above a certain threshold if $|V(x)| \leq C/(x+1).$ In \cite{KLS}, it was shown 
that any imbedded eigenvalues $\lambda_n$ that can occur for $|V(x)| \leq C/(x+1)$ satisfy $\sum_n \lambda_n <  \infty.$ 
This shows that a construction carried out in this paper is not  possible for $V$ decaying at a Coulomb rate. We conjecture that 
singular continuous spectrum cannot occur for such $V.$ \\


{\bf Acknowledgment. } The author is grateful to Michael Christ and 
Barry Simon for useful and stimulating discussions. 
The author also wishes to express his gratitude 
to Francois Germinet and Stephan de Bi\`evre at Universit\'e de Lille I, 
where part of this work was carried out.    


\end{document}